\documentclass[10pt]{article}

\newtheorem{prop}{Proposition}
\newtheorem{theoreme}{Theorem}
\newtheorem{lemme}{Lemma}
\newcommand{\eqref}[1]{\mbox{{(\ref{#1})}}} 
\newcommand{\norme}[2]{ | \, #1 \, |_{#2}}
\newcommand{\gnorme}[2]{ | \,  #1 \, |_{ w , #2}}
\newcommand{\hnorme}[2]{| \,  #1 \, |_{H^{#2}} }
\newcommand{\hhnorme}[2]{| \,  #1 \, |_{\dot{H}^{#2}} }
\newcommand{\R}[1]{\mathbf{R}^{#1}}
\newcommand{\N}[1]{\mathbf{N}^{#1}}
\newcommand{\IInt}[1]{\int \! \! \! \! \int_{#1}}
\newenvironment{proof}{$ $ \\ \textbf{Proof.}} { \hspace{\stretch{1}}
 \rule{1ex}{1ex} \\ $ $ \\}

\title{Asymptotic stability of Oseen vortices for a density-dependent incompressible viscous fluid}

\author{{\bf L.Miguel Rodrigues}\\
Institut Fourier \\ U.J.F. Grenoble 1, C.N.R.S.\\
100 rue des maths\\38402 Saint Martin d'H\`eres, France}

\date{}

\begin{document}

\maketitle

{\hspace{\stretch{1}}} {\bf Abstract} {\hspace{\stretch{1}}}

In the analysis of the long-time behaviour of two-dimensional
incompressible viscous fluids, Oseen vortices play a major role as
attractors of \emph{any} homogeneous solution with integrable initial
vorticity \cite{Gallay_W-global_stability}. As a first step in the
study of the \emph{density-dependent} case, the present paper
establishes the asymptotic stability of Oseen vortices for slighty
inhomogeneous fluids with respect to localized regular perturbations.\\



\section*{Introduction}

In this paper we consider the motion of a weakly inhomogeneous
incompressible viscous fluid in the two-dimensional Euclidean space. We
can modelize the fluid by $(\rho,u)$, where $\rho=\rho(t,x) \in \R +$
is the density field and $u=u(t,x) \in \R 2$ the velocity field. The
evolution we consider here is governed by the density-dependent
incompressible Navier-Stokes equations: 
\begin{equation} \label{inhomogeneNS-velocity}
\left\{
\begin{array}{rcl}
\partial _{t} \, \rho + (u \cdot \nabla) \, \rho & = & 0 \\
\partial _{t} \,  u + (u \cdot \nabla) \, u & = & 
\frac{1}{\rho} \, ( \triangle u - \nabla p ) \\
\textnormal {div} \, u & = & 0
\end{array}
\right.
\end{equation}
where $p=p(t,x) \in \R{}$ is the pressure field, which is determined
(up to a constant) by the preservation of incompressibility which yields the elliptic equation:
\begin{equation} \label{pression}
\textnormal{div}\,\big(\,\frac{1}{\rho}\,\nabla p\,\big) 
=\textnormal{div}\,\big(\,\frac{1}{\rho}\,\triangle\,u-(u \cdot \nabla)\,u\,\big) \ .
\end{equation}

Alternatively, we can represent the fluid motion using the vorticity
field $\omega = \textnormal{curl} \ u \in \R{}$ rather than the
velocity. Note that, in the two-dimensional context, $\textnormal{curl}
(f_1,f_2)$ stands for $\partial_1 f_2 - \partial_2 f_1$. Therefore the
evolution equations for $(\rho,\omega)$ become
\begin{equation} \label{inhomogeneNS-vorticity}
\left\{ \begin{array}{rcl}
\partial _{t} \, \rho + (u \cdot \nabla) \, \rho & = & 0 \\
\partial _{t} \, \omega + (u \cdot \nabla) \, \omega & = & 
\textnormal{div}\,\big( \, \frac{1}{\rho} \,
( \nabla \omega + \nabla^{\perp} p) \, \big)
\end{array} \right.
\end{equation}
where $p$ is again determined by \eqref{pression}, and $u$ is 
recovered from $\omega$ \emph{via} the Biot-Savart law:
\begin{equation} \label{BS}
u(x) = \frac{1}{2 \pi} \int_{\R{2}}
\frac{(x-y)^{\perp}}{| x-y | ^{2}} \, \omega (y) \, dy
\end{equation}
for $x \in \R2$, with $(z_1, z_2) ^{\perp} = (-z_2, z_1)$. We also
note $u = K_{BS} \star \omega$, where $K_{BS}$ is the Biot-Savart kernel:
$K_{BS} (x) =  \frac{1}{2 \pi} \frac{x^{\perp}}{| x | ^{2}} $. For
simplicity, throughout the present paper, the viscosity of the fluid
has been rescaled to one.  

We refer to the monograph \cite{Lions1} for a general presentation of
the available mathematical results on incompressible Navier-Stokes
equations. We also mention the work of B.Desjardins on the global
existence of weak solutions
\cite{Desjardins-transport_NS,Desjardins-NS_inhomogene}, and, closer
to the spirit of the present paper, the work of R.Danchin on
well-posedness in Besov spaces \cite{Danchin-NS_inhomogene}. Let us
emphasize that they both work with the velocity formulation
\eqref{inhomogeneNS-velocity} and do not assume the density $\rho$ to
be bounded away from zero. In more physical terms, they allow for
regions of (almost complete) vacuum, which create technical difficulties.

In contrast, not only shall we not allow the density to be close to
zero but we shall only consider \emph{weakly} inhomogeneous fluids,
namely we shall assume that the density $\rho$ is close to a positive
constant which, without loss of generality, we take equal to
one. Remark that if the initial density is contant in space
\emph{i.e.} homogeneous, then the density remains equal to this
constant for all subsequent times. Therefore, in such a case, system
\eqref{inhomogeneNS-velocity} reduces to the usual incompressible
Navier-Stokes equations. Moreover, since $\textnormal
{div}(\nabla^{\perp} p) =0$, the pressure term disappears from system
\eqref{inhomogeneNS-vorticity} which thus reduces to 
\begin{equation} \label{NS-vorticity}
\partial _{t} \, \omega + (u \cdot \nabla) \, \omega = \triangle \, \omega \ .
\end{equation}

Again, a wealth of information on the Cauchy problem for the
\emph{homogeneous} incompressible Navier-Stokes equations can be found
in \cite{Lions1} or \cite{Lemarie}. Concerning the long-time behaviour
of the solutions of the vorticity equation \eqref{NS-vorticity}, the
work of Th.Gallay and C.E.Wayne has revealed the important role played
by a family of explicit self-similar solutions, \emph{Oseen
  vortices}, given by $\rho \equiv 1$, $u = \alpha \, u^{G}$ and
$\omega = \alpha \, \omega^{G}$, where 
\begin{displaymath}
\omega^G (t,x)  =  \frac{1}{t} \, G \big( \frac{x}{\sqrt{t}} \big) 
\ , \quad
u^G (t,x) =  \frac{1}{\sqrt{t}} \, v^{G} \big( \frac{x}{\sqrt{t}}\big)
\end{displaymath}
with
\begin{displaymath}
\qquad
G (\xi ) = \frac{1}{4 \pi} \, e^{- | \xi |^{2} / 4} 
\  , \quad
v^{G} (\xi ) = \frac{1}{2 \pi} \, \frac{\xi^{\perp}}{|\xi|^{2}}
\, \big( 1 - e^{- | \xi |^{2} / 4} \big)
\end{displaymath}
and $\alpha \in \R{}$ parametrizes this family of solutions. For an
Oseen vortex $(1,\alpha\,\omega^G)$, $|\alpha|$ is actually its Reynolds number.
If the initial vorticity $\omega_0$ is integrable, it is proved in
\cite{Gallay_W-global_stability} that the corresponding solution of
\eqref{NS-vorticity} converges to $\alpha \, \omega^G$ in $L^1$-norm
as $t \to \infty$, where $\alpha :=  \int_{\R 2} \omega_0$. Moreover,
it was shown in \cite{G_Gallay-uniqueness,G_Gallay_Lions-uniqueness}
that $\alpha \, \omega^G$ is the unique solution of the vorticity
equation \eqref{NS-vorticity} with initial data $\alpha \, \delta_0$.

Even though the homogeneous incompressible Naviers-Stokes equations
provide a good model in many situations, all real fluids are, at least
slightly, inhomogeneous and it is therefore important and relevant,
both from a practical and theorical point of view, to investigate
whether the predictions of the homogeneous model are meaningful for
the \emph{density-dependent} model, especially in the weakly
inhomogeneous regime. The goal of this paper is to address this
question in the particular case of Oseen vortices. These explicit
solutions persist in the density-dependent case if we
assume $\rho \equiv 1$, and it is therefore natural to ask whether
they play there the same role as in the homogeneous case. While the
general answer to this question is unknown, we treat here one
important aspect: are Oseen vortices stable solutions for the
density-dependent incompressible Navier-Stokes equations ? In other
words, does the theory predicts that these self-similar solutions
may be observed ?

Before stating what we mean exactly by stability, let us recall an
important property of the Navier-Stokes equations: \emph{scaling
  invariance}. For any $\lambda >0$, if $(\rho(t,x),\omega(t,x))$ is a
solution of \eqref{inhomogeneNS-vorticity}, with corresponding
velocity $u(t,x)$ and pressure $p(t,x)$, so is 
\begin{displaymath}
D_{\lambda} (\rho,\omega) = (\rho \, (\lambda^2 t,\lambda x),
\lambda^2 \, \omega \, (\lambda^2 t,\lambda x)) \ , 
\end{displaymath}
with velocity $\lambda \, u(\lambda^2 t,\lambda x)$ and pressure
$\lambda^2 \, p(\lambda^2 t,\lambda x)$. Moreover, Oseen vortices are
\emph{self-similar}, in the sense that $D_{\lambda} (1,\alpha \,
\omega^G) = (1,\alpha \, \omega^G)$, for any $\alpha \in \R{}$ and any
$\lambda >0$. To study these solutions, it is therefore more
convenient to rewrite \eqref{inhomogeneNS-vorticity} in self-similar
variables. Following \cite{Gallay_W-invariant_manifold} , we set
\begin{equation} \label{selfsimilar}
(\tau,\xi) := (\ln t, \frac{x}{\sqrt{t}}) \ .
\end{equation}
Motivated by scaling invariance, we will work with new quantities
$(r,w,v,\Pi)$ related to the former by
\begin{equation} \label{selfsimilar2}
\begin{array}{rclcrclc}
r(\tau, \xi)   &=& \rho \, (e^{\tau},e^{\frac{\tau}{2}} \xi) &,&
v(\tau, \xi)   &=& 
e^{\frac{\tau}{2}} \, u \, (e^{\tau},e^{\frac{\tau}{2}} \xi) &,\\ 
w(\tau, \xi)   &=& e^{\tau} \, \omega \, (e^{\tau},e^{\frac{\tau}{2}} \xi) &,&
\Pi(\tau, \xi) &=& e^{\tau} \, p \, (e^{\tau},e^{\frac{\tau}{2}} \xi) &.
\end{array}
\end{equation}
Then the corresponding evolution equations for $(r,w)$ are
\begin{equation} \label{inhomogeneNS-selfsimilar}
\left\{ \begin{array}{rcl}
\partial_{\tau} \, r + 
\big( (v - \frac{1}{2} \, \xi) \cdot \nabla \big) \, r
& = & 0 \\
\partial_{\tau} \, w + 
\big( (v - \frac{1}{2} \, \xi) \cdot \nabla \big) \, w
- w & = &
\textnormal {div} \, \big( \, \frac 1r \,
( \nabla \omega + \nabla^{\perp} \Pi) \, \big) \ 
\end{array}
\right.
\end{equation}
where again $v$ is obtained from $w$ by the Biot-Savart law and
$\nabla \Pi$ by solving
\begin{equation} \label{pression-selfsimilar}
\textnormal {div} \, \big( \, \frac 1r  \, \nabla \Pi \, \big) 
= \textnormal {div} \, \big( \, \frac 1r \, \triangle \, v - (v \cdot
\nabla) \, v \big) \ .
\end{equation}
By constuction, Oseen vortices correspond to stationary solutions of
\eqref{inhomogeneNS-selfsimilar}, namely $(1, \alpha \, G)$, where
$\alpha$ is any real number. We fix initial data for the original
equations at $t=1$ rather than at $t=0$, hence at $\tau=0$ for the
new equations.

In order to state our main result, we now write down the evolution
equations for a perturbation of an Oseen vortex. Let $\alpha \in
\R{}$. We work with new quantities $(b,\widetilde{w})$ related to the
former by $b = \frac 1r -1$ and $\widetilde{w} = \alpha \, G - w$. The
evolution equations for $(b,\widetilde{w})$ are
\begin{equation} \label{complete}
\left\{ \begin{array}{rcl}
\partial_{\tau} \, b + 
\big( (v - \frac 12 \, \xi) \cdot \nabla \big) \, b
& = & 0 \\
\partial_{\tau} \, \widetilde{w} 
- ( \mathcal{L} - \alpha \, \Lambda ) \, \widetilde{w} 
+ \big( \widetilde{v} \cdot \nabla \big) \, \widetilde{w} 
& = & \textnormal{div} \, \big( \, b \, ( \nabla w
+ \nabla^{\perp} \Pi ) \, \big )
\end{array}
\right.
\end{equation}
where 
\begin{equation}\label{LandLambda}
\begin{array}{rcl}
\mathcal{L} f & = & \triangle f + \frac 12 \, \xi \cdot \nabla f + f 
\ ,\\
\Lambda f & = & v^G \cdot \nabla f + (K_{BS} \star f) \cdot \nabla G \ .
\end{array}
\end{equation}
Here $\widetilde{v}$ is obtained from $\widetilde{w}$ by the Biot-Savart
law, $w$ and $v$ are recovered by
\begin{equation} \label{perturbation}
w = \alpha \, G + \widetilde{w} \ , \quad
v = \alpha \, v^G + \widetilde{v} \ ,
\end{equation}
and $\nabla \Pi$ is obtained by solving
\begin{equation} \label{pression-perturbation}
\textnormal{div} \, \big( \, (1 + b) \, \nabla \Pi  \, \big)  
 = \textnormal{div} \, \big( \, (1 + b) \, \triangle v 
- (v \cdot \nabla) \, v  \,\big) \ .
\end{equation}
Let us point out that, thanks to the linearity of the Biot-Savart law,
we do have $v=K_{BS} \star w$ and $\widetilde{v}=K_{BS} \star \widetilde{w}$. 

Note that $\mathcal{L}$ is the usual Fokker-Planck type operator which
generates  the evolution corresponding to the heat equation in
self-similar variables. On the other hand, $\Lambda$ is the non-local
first-order operator resulting from the linearization around $w=G$ of
the transport term of \eqref{inhomogeneNS-selfsimilar}. More
precisely, if $\widetilde{v}=K_{BS} \star \widetilde{w}$ and $(v,w)$
satisfies \eqref{perturbation}, then $v \cdot \nabla w= \alpha \,
\Lambda \, \widetilde{w} + \widetilde{v} \cdot \nabla \widetilde{w}$,
since $v^G \perp \nabla G$.

Before stating our result of stability, we finally introduce the
function spaces and norms we will encounter throughout the present
paper. For $1 \leq p \leq \infty$, we denote by $L^p (\R 2)$ the usual
$L^p$-space and by $\norme{f}{p}$ the usual $L^p$-norm of a function
$f$. Similarly, for $s \in \R{}$, we denote by $H^s (\R 2)$ the usual
Sobolev space and by $\hnorme{f}{s}$ the corresponding norm. Sobolev norms
will be convenient to specify the regularity of our perturbations, but
we also need weighted norms to describe the localization
properties. Indeed, even in the homogeneous case, we know that it is
impossible to get a convergence rate in time if we do not assume that
the perturbations decay sufficiently fast at infinity in space (see
\cite{Gallay_W-global_stability}). Instead of using polynomial weights
as in \cite{Gallay_W-global_stability}, we choose here the Gaussian
weight $G^{-\frac 12}$, which is naturally related to the Oseen
vortices and has several technical advantages. For instance, on the
Hilbert space $L_w^2(\R2):=\{ f\ |\ G^{-\frac 12}\,f \in L^2(\R2) \}$,
the linear operators $\mathcal{L}$ and $\Lambda$ become respectively
symmetric and skew-symmetric. Consequently, for any $1\leq p \leq
+\infty$, we shall use the weighted $L^p$-space defined as follows: 
\begin{equation} \label{GLp}
L_w^p (\R2):=\{ \ f\ |\ G^{-\frac 12}\,f \in L^p(\R2) \ \}
\end{equation}
with the corresponding norms $\gnorme{f}{p}:=\norme{G^{-\frac 12}\,f}{p}$.

We are now able to state the main result of this paper:

\begin{theoreme} \label{main}
Let $\alpha \in \R{}$, $0<s<1$ and $0<\gamma<\frac 12$. 
There exist $\varepsilon_0 >0$ and $K>0$ such that for any
$0<\varepsilon<\varepsilon_0$ if
\begin{enumerate}
\item $b_0 \in H^{s+2} (\R 2)$ is such 
that $G^{-\frac 12} b_0$ belongs to $L^2(\R 2) \cap L^{\infty}(\R 2)$
with $\gnorme{b_0}{2} \leq \varepsilon$, 
$\gnorme{b_0}{\infty} \leq \varepsilon$, 
and $G^{- \frac 12} \nabla b_0$ to $L^q (\R 2)$ 
for some $q > \max(4,\frac2s)$
\item $\widetilde{w}_0 \in H^s (\R 2)$ is such that 
$G^{-\frac 12} \widetilde{w}_0$ belongs to $L^{2} (\R 2)$ with 
$\gnorme{\widetilde{w}_0}{2} \leq \varepsilon$, and
$\int_{\R2} \widetilde{w}_0 = 0$
\end{enumerate}
then there exists a unique solution $(b,\widetilde{w})$ of \eqref{complete}
with initial data $(b_0,\widetilde{w}_0)$ such that
\begin{enumerate}
\item $b \in L^{\infty}_{loc} (\R +;H^{s+2}(\R 2))$, 
\item $G^{-\frac12} b \in L^{\infty}(\R+;L^2(\R 2) \cap L^{\infty}(\R 2))$,
$G^{-\frac12} \nabla b \in L^{\infty}_{loc} (\R +;L^q (\R 2))$
\item $\widetilde{w} \in L^{\infty}_{loc}(\R +;H^{s}(\R 2)) \cap 
L^{2}_{loc}(\R +;H^{s+1}(\R 2))$,
\item $G^{-\frac12}\widetilde{w} \in  L^{\infty}(\R +;L^2(\R 2)) \cap 
L^2(\R +;L^2(\R 2))$,
$G^{-\frac12} \nabla \widetilde{w} \in  L^2(\R +;L^2(\R 2))$,
$G^{-\frac12} | \xi | \, \widetilde{w} \in  L^2(\R +;L^2(\R 2))$.
\end{enumerate}
Moreover this solution satisfies $\gnorme{\widetilde{w} (\tau)}{2} \leq K
\, \varepsilon \, e^{- \gamma \, \tau}$, for any $\tau>0$.
\end{theoreme}

Let us make some comments on this result.

\begin{enumerate}
\item We treat only small regular and localized perturbations. Note also that one can not allow for perturbations with non-zero mean vorticity. Indeed, since $\int_{\R 2} \widetilde{w}$ is preserved under the evolution, in order to have convergence in $L^1(\R2)$, we must have 
$\alpha=\int_{\R 2} w_0$. Consequently, henceforth, we will scarcely repeat but always assume $\int_{\R2} \widetilde{w}_0 = 0$.

\item In contrast, one important aspect of our result is that we make no smallness assumption on the size of the Oseen vortices, so we do treat high Reynolds numbers $|\alpha|$. Yet note again that the present paper is far from telling us whether Oseen vortices attract \emph{all} weakly inhomogeneous fluids.

\item Now, to make our result more concrete, let us give some immediate consequences in the original variables. Under the asumptions of {\bf Theorem~\ref{main}}, if $(\rho,\omega)$ is the solution of \eqref{inhomogeneNS-vorticity} defined by \eqref{selfsimilar2} with $w=\alpha\,G+\widetilde{w}$ and $r=\frac{1}{1+b}$, then the vorticity $\omega(t,x)$ satisfies
\begin{displaymath}
t^{\frac12} \, \norme{e^{\frac{|\xi|^2}{8 t}} \, 
(\omega(t)-\alpha\,\omega^G(t))}{2}\, \leq \,\frac{C}{t^{\gamma}} \ ,
\quad  t\geq 1 \ .
\end{displaymath}
Moreover this implies
\begin{displaymath}
t^{1-\frac1p} \, \norme{\omega(t)-\alpha \, \omega^G(t)}{p} \,
\leq \, \frac{C_p}{t^{\gamma}}, \quad
t^{\frac12-\frac1q} \, \norme{u(t)-\alpha \, u^G(t)}{q} \,
\leq \, \frac{C_q}{t^{\gamma}}
\end{displaymath}
for any $1 \leq p \leq 2$, \, $2< q <+\infty$, \, $t \geq 1$. Note that $t^{1-\frac1p}$ and $t^{\frac12-\frac1q}$ are the rates of convergence corresponding to Oseen vortices.

\item Concerning the density, we can not expect the perturbation to decay in time in the whole space, since we deal with an incompressible fluid. Indeed, since the density satisfies a transport equation by a divergence-free velocity field, formally it remains constant in law. Thus any $L^p$-norm of a function of $\rho$ remains constant in time.
\end{enumerate}

To establish {\bf Theorem~\ref{main}}, following \cite{Danchin-NS_inhomogene}, we build a sequence of solutions of a linearization of
\eqref{complete}. In order to show the convergence of this sequence,
we will use local-in-time estimates in Sobolev norms for solutions of
this linearization. To be more precise, in a \emph{first step}, we
establish global estimates for solutions of the linearized system,
which control the density in weighted $L^p$ spaces and the vorticity
in $H^1$. In a \emph{second step}, we use the previous estimates to
prove first global estimates of the density in $H^{2-\varepsilon}$
(with a loss of regularity), then local estimates of the vorticity in
$H^{s+1}$, and finally local estimates of the density in
$H^{s+2}$. Then we use these estimates to establish the existence and
uniqueness parts of the theorem. Finally, we prove the asymptotic part of
the theorem. 

The plan of the paper does not follow the former steps, but we will
encounter the structure of the proof in the various sections. In a
preliminary section, we collect estimates on the Biot-Savart kernel,
thus on the velocity in terms of the vorticty, and on the pressure. In
the second section, we gather estimates on the linearized density
equation. In the third one, we establish the estimates on the
linearized vorticity equation. At last, we state and
prove our main results. 

In what follows, the original (and physical) time will never appear
again, so for notational convenience henceforth we allowed us to use the letter $t$ to denote the rescaled time $\tau$.

\section{Preliminaries}

If $f$ is integrable over $\R2$, we define its Fourier transforms to be the function $\hat{f}$ defined for any $\eta \in \R 2$ by
\begin{displaymath}
\hat{f} (\eta) = \int_{\R 2} \, f(\zeta) \, e^{i \eta \cdot \zeta} \,
d\zeta \ .
\end{displaymath}

Concerning function spaces, we will also need the following
convention. For any $\sigma \in \R{}$, we denote by $\dot{H}^{\sigma}(\R2)$ the usual homogeneous Sobolev space on $\R 2$ equiped with $\hhnorme{f}{\sigma}:=\norme{I^{\sigma} f}{2}$, where $I:=(-\triangle)^{\frac12}$.

\subsection{Biot-Savart kernel}

The goal of this subsection is to enable us to estimate the velocity in terms of the vorticity. We collect here some estimates on $v$ in terms of $w$ when $v$ is obtained from $w$ by the Biot-Savart law, thus on the velocity in terms of the vorticity. Recall that, in this case, for almost every $\xi \in \R 2$,
\begin{equation} \label{BiotSavart}
 v (\xi) = \frac{1}{2 \pi} \int_{\R{2}}
\frac{( \xi- \eta )^{\perp}}{| \xi - \eta | ^{2}} \, w (\eta) \, d\eta 
\end{equation}
where $(x_{1} , x_{2} ) ^{\perp} = (-x_{2} , x_{1})$, that is $v = K_{BS} \star w$, where $K_{BS}$ is the Biot-Savart kernel 
$K_{BS} (x) = \frac{1}{2 \pi} \frac{x^{\perp}}{| x | ^{2}} $. Also note that in terms of Fourier transform, \eqref{BiotSavart} becomes \begin{equation} \label{BiotSavartFourier} 
\hat{v}(\eta) =\frac{i\,\eta^{\perp}}{|\eta|^2}\, \hat{w}(\eta) \, .
\end{equation}
Most of these estimates are already known, but for the sake of completeness
we try to give proofs rather than references.

The following proposition gathers estimates in $L^p$-spaces.

\begin{prop} $ $
\begin{enumerate}
\item Let $1 < p < 2 < q$ be such that $1 + \frac 1q= \frac 12 + \frac 1p$.\\
There exists $C > 0$ such that, if $w$ belongs to $L^{p} (\R{2})$, 
then \eqref{BiotSavart} defines a $v$ in $L^{q} (\R{2})$ and
\begin{equation} \label{BiotSavartLp} 
\norme{v}{q}  \leq C \, \norme{w}{p} \ .
\end{equation}

\item Let $1 \leq p < 2 < q \leq + \infty$ and $ 0 < \theta < 1$ be  such that $\frac{\theta}{p} + \frac{1-\theta}{q} = \frac 12$.\\
There exists $C > 0$ such that, if $w$ belongs to 
$L^{p} (\R{2}) \cap L^{q} (\R{2})$, 
then \eqref{BiotSavart} defines a $v$ in $L^{\infty} (\R{2})$ and
\begin{equation} \label{BiotSavartL8} 
\norme{v}{\infty} \leq C\,\norme{w}{p}^{\theta}\,\norme{w}{q}^{1-\theta}
\ .
\end{equation}

\item Let $p>1$.
There exists $C > 0$ such that, if $w$ belongs to $L^{p} (\R{2})$ and
$v$ is defined from $w$ by \eqref{BiotSavart}, then $\nabla v$ belongs
to $L^{p} (\R{2})$ and
\begin{equation} \label{BiotSavartGrad} 
\norme{\nabla v}{p}  \leq C \, \norme{w}{p} \ .
\end{equation}

\end{enumerate}
In addition, in all cases, we have $\textnormal{div} \, v = 0$ and 
$ \textnormal{curl}  \, v = w$.
\end{prop}

\begin{proof}
Part 1 follows from a Young-like inequality called Hardy-Littlewood-Sobolev inequality (see for instance {\bf Theorem~V.1} in \cite{Stein}). Indeed, $K_{BS}$ is weakly $L^2$ but not square integrable.

Part 2 is trivial when $w \equiv 0$. If not, we remark that from H\"older's inequalities, we obtain 
\begin{displaymath}
\begin{array}{rcl}
|v(\xi)| &\leq& 
\frac{1}{2 \pi}\, \int_{\{|\eta|\leq R\}} |w(\xi-\eta)|\,\frac{1}{|\eta|} \, d\eta
+\frac{1}{2 \pi}\, \int_{\{|\eta|\geq R\}} |w(\xi-\eta)|\,\frac{1}{|\eta|} \, d\eta\\
 &\leq& C\ \norme{w}{q}\ R^{1-\frac2q} 
      + C\ \norme{w}{p}\ \frac{1}{R^{\frac2p -1}} \ ,
\end{array}
\end{displaymath}
for almost every $\xi \in \R 2$ and any $R >0$. Aiming to optimize this inequality, we choose $R=(\frac{\norme{w}{p}}{\norme{w}{q}})^{\beta}$ with 
$\beta = \frac{1-\theta}{\frac2p - 1} = \frac{\theta}{1-\frac2q}$ and derive \eqref{BiotSavartL8}.

Part 3 holds since derivating \eqref{BiotSavart} yields that $\nabla
v$ is obtained from $w$ by a singular integral kernel of
Calder\'on-Zygmund type (see {\bf Theorem II.3} in \cite{Stein}).
\end{proof}

The following proposition gathers estimates in Sobolev spaces.

\begin{prop} $ $
\begin{enumerate}
\item Let $s \in \R{}$. There exists $C>0$ such that, if $w$ belongs to 
$\dot{H}^{s-1} (\R2)$, then, if $v$ is defined by
\eqref{BiotSavart}, $v$ belongs to $\dot{H}^s (\R2)$ and
\begin{equation} \label{BiotSavartHsGrand} 
\hhnorme{v}{s}  \leq C \, \hhnorme{w}{s-1} \  .
\end{equation}

\item Let $0<s<1$. \\
There exists $C > 0$ such that, if $(1+|\xi|)\,w$ belongs to 
$L^2(\R2)$, then, if $v$ is defined by \eqref{BiotSavart}, $v$
belongs to $\dot{H}^s(\R2)$ and
\begin{equation} \label{BiotSavartHsPetit} 
\hhnorme{v}{s} \leq 
C\,\norme{(1+|\xi|)\,w}{2} \ . 
\end{equation}

\end{enumerate}
\end{prop}

\begin{proof}
Part 1 is a direct consequence of the Fourier formulation
\eqref{BiotSavartFourier}.

Part 2 is thus reduced to estimate $\hhnorme{w}{s-1}$. If $0<s<1$, we note that
\begin{displaymath}
\hhnorme{w}{s-1}^2 =
C\,\int_{\R2} \,\frac{|\hat{w}(\xi)|^2}{|\xi|^{2(1-s)}}\, d\xi
\ \leq \ 
C\,\int_{|\,\xi\,|\leq 1} \,\frac{|\hat{w}(\xi)|^2}{|\xi|^{2 (1-s)}}\, d\xi
+C\,\int_{|\,\xi\,|\geq 1} \,|\hat{w}(\xi)|^2\, d\xi \ .
\end{displaymath}
The second term of the right member is dominated by $\norme{w}{2}^2$. Choosing $p$ such that $p>\frac 2s$ and applying first H\"older's inequalities then Sobolev's embeddings, we obtain for the first term
\begin{displaymath}
\int_{|\,\xi\,| \leq 1} \,\frac{|\hat{w}(\xi)|^2}{|\xi|^{2 (1-s)}}\, d\xi
\ \leq \ C\,\norme{\hat{w}}{p}^2
\ \leq \ C\,\hnorme{\hat{w}}{1}^2 \ .
\end{displaymath}
Finally, gathering every piece yields
\begin{displaymath}
\hhnorme{w}{s-1}^2  
\ \leq \ C\,\hnorme{\hat{w}}{1}^2+C\,\norme{w}{2}^2 
\ \leq \ C\,\norme{(1+|\xi|)\,w}{2}^2 \ .
\end{displaymath}
This concludes the proof.
\end{proof}

\subsection{Pressure estimates}

Keeping in mind equation \eqref{pression-perturbation}, we gather
some estimates for a solution $\Pi$ of the following equation:
\begin{equation} \label{pression-type}
\textnormal{div} \, \big( (1 + b) \nabla \Pi \big)   =
\textnormal{div}\, F \ .
\end{equation}

We begin with estimates in $L^p$-spaces.

\begin{prop} $ $
\begin{enumerate}
\item Let $p >1$. \\
There exists $C > 0$ and $\kappa > 0$ such that if $F$ belongs to
$L^{p} (\R{2})$ and $b$ to $L^{\infty} (\R{2})$ with $\kappa \,
\norme{b}{\infty} < 1$, then \eqref{pression-type} has a unique solution
$\Pi$ (up to a constant) such that $\nabla \Pi$ belongs to $L^{p}(\R{2})$,
and 
\begin{equation} \label{pressionLp}
\norme{\nabla \Pi}{p} 
\leq \frac{C}{1-\kappa \, \norme{b}{\infty}} \, \norme{F}{p} \ .
\end{equation}

\item Let $1 < p, q, r$ be such that 
$\frac 1r = \frac 1p + \frac 1q$.\\
There exists $C > 0$ and $\kappa > 0$ such that if $F$ belongs to
$L^{q} (\R{2})$ and, for $i=1,2$, $b_{i}$ belongs to $L^{\infty}
(\R{2}) \cap L^{p} (\R{2})$  with $\kappa \, \norme{b_{i}}{\infty} <
1$ and $\Pi_{i}$ solves 
\begin{displaymath}
\textnormal{div} \, \big( (1 + b_{i}) \nabla \Pi _{i} \big)   =
\textnormal{div} \, F \ ,
\end{displaymath}
then $\nabla \, (\Pi_{2} - \Pi_{1})$ belongs to $L^{r} (\R{2})$ and 
\begin{equation} \label{pressiondelta}
\norme{\nabla \, (\Pi_{2} - \Pi_{1})}{r} 
\leq \frac{C}{(1-\kappa \, \norme{b}{\infty})^{2}} \, 
\norme{b_{2} - b_{1}}{p} \, \norme{F}{q} \ .
\end{equation}
\end{enumerate}
\end{prop}

\begin{proof}
\begin{enumerate}
\item  We want to obtain $\nabla \Pi$, the solution of 
\eqref{pression-type}, in terms of $F$ as a perturbation of Leray projectors. Let $\mathbf{P}$ be the Leray projector, that is the projector on divergence-free vector fields along gradients, and let $\mathbf{Q}= \mathbf{I} - \mathbf{P}$. Remark that $\mathbf{Q} \, F$ gives the solution $\nabla \Pi$ of \eqref{pression-type} when $b \equiv 0$. Now we can rewrite \eqref{pression-type} as
\begin{displaymath}
\textnormal{div} \, \nabla \Pi = 
\textnormal{div} \, ( \, F - b \, \nabla \Pi \, )
\end{displaymath}
then $\nabla \Pi=\mathbf{Q}\,(F-b\,\nabla \Pi)$ thus
$(\mathbf{I}+\mathbf{Q}\,b)\,\nabla \Pi=\mathbf{Q}\,F$. Since $\mathbf{Q}$ is continuous on $L^p$, there exists $\kappa>0$ such that
\begin{displaymath}
\norme{\mathbf{Q}\,b\,f}{p}\,\leq\,\kappa\,\norme{b}{\infty}\,\norme{f}{p}
\end{displaymath}
for $f \in L^p(\R2)$. Thus, if $\kappa\,\norme{b}{\infty}<1$, \ $\mathbf{I}+\mathbf{Q}\,b$\ is invertible on $L^p$, and
\begin{equation} \label{Leray}
\nabla \Pi=(\mathbf{I}+\mathbf{Q}\,b)^{-1}\,\mathbf{Q}\,F
\end{equation}
gives the unique solution, with the expected bound.

\item Reminding \eqref{Leray} we write
\begin{displaymath}
\begin{array}{rcl}
\nabla \Pi_1 &=& (\mathbf{I}+\mathbf{Q}\,b_2)^{-1}\,
(\mathbf{I}+\mathbf{Q}\,b_2)\,(\mathbf{I}+\mathbf{Q}\,b_1)^{-1}\,
\mathbf{Q}\,F \\
\nabla \Pi_2 &=& (\mathbf{I}+\mathbf{Q}\,b_2)^{-1}\,
(\mathbf{I}+\mathbf{Q}\,b_1)\,(\mathbf{I}+\mathbf{Q}\,b_1)^{-1}\,
\mathbf{Q}\,F 
\end{array}
\end{displaymath}
then subtracting and factorizing yields
\begin{displaymath}
\nabla (\Pi_2-\Pi_1)= (\mathbf{I}+\mathbf{Q}\,b_2)^{-1}\,\mathbf{Q}\,
(b_1-b_2)\,(\mathbf{I}+\mathbf{Q}\,b_1)^{-1}\,\mathbf{Q}\,F \ .
\end{displaymath}
Now the continuity of the operator $\mathbf{Q}$ on $L^r$ reduces 
\eqref{pressiondelta} to an estimate on 
$\norme{(b_1-b_2)\,(\mathbf{I}+\mathbf{Q}\,b_1)^{-1}\,\mathbf{Q}\,F}{r}$. At last applying first H\"older's inequalities then the continuity on 
$L^q$ concludes the proof.
\end{enumerate}
\end{proof}

In order to estimate solutions of \eqref{pression-type} in Sobolev spaces,
we first state some useful commutator estimates of Kato-Ponce
type (see \cite{KatoPonce}). Let us recall that 
$I=(-\triangle)^{\frac 12}$.

\begin{lemme} Let $0<s<1$ and $\sigma >1$.
\begin{enumerate}
\item There exists $C>0$ such that if $I^s f$ belongs to $L^2(\R2)$ and 
$g$ to $H^{\sigma}(\R2)$, then $I^s(fg)-f\,I^s g$ belongs to $L^2(\R2)$ and 
\begin{equation} \label{commutateurHs1}
\norme{I^s(fg)-f\,I^s g}{2} 
\leq C\,\norme{I^s f}{2}\,\hnorme{g}{\sigma}.
\end{equation}

\item There exists $C>0$ such that if $I^s f$ belongs to $H^{\sigma}(\R2)$ and $g$ to $L^2(\R2)$, then $I^s(fg)-f\,I^s g$ belongs to $L^2(\R2)$ and 
\begin{equation} \label{commutateurHs2}
\norme{I^s(fg)-f\,I^s g}{2} 
\leq C\,\hnorme{I^s f}{\sigma}\,\norme{g}{2}.
\end{equation}
\end{enumerate}
\end{lemme}

\begin{proof}
Let us first note that there exists $C >0$ such that
\begin{displaymath}
\norme{\hat{h}}{1} \, \leq \, C \, \hnorme{h}{\sigma} \ ,
\end{displaymath}
for any $h$ in $H^{\sigma} (\R 2)$. This comes applying H\"older's
inequalities to 
\begin{displaymath}
\int_{\R 2} |\hat{h}(\eta)|\,d\eta = 
\int_{\R 2} \, \frac{1}{(1+|\eta|^2)^{\frac{\sigma}{2}}} \,
(1+|\eta|^2)^{\frac{\sigma}{2}} \, |\hat{h}(\eta)|\,d\eta \ .
\end{displaymath}

Therefore in order to prove the lemma it is sufficient to establish
\begin{eqnarray}
\norme{I^s(fg)-f\,I^s g}{2} 
&\leq&C\,\norme{I^s f}{2}\,\norme{\hat{g}}{1}\ , \label{commHs1} \\
\norme{I^s(fg)-f\,I^s g}{2} 
&\leq&C\,\norme{\widehat{I^s f}}{1}\,\norme{g}{2}\ . \label{commHs2}
\end{eqnarray}

Now set $h=I^s(fg)-f\,I^s g$. We have
\begin{displaymath}
\hat{h}(\eta) = \frac{1}{(2\pi)^2} \int_{\R 2}\,(|\eta|^s-|\eta-\zeta|^s)\,
\hat{f}(\zeta)\,\hat{g}(\eta-\zeta)\,d\zeta
\end{displaymath}
for almost every $\eta \in \R 2$. Thanks to the following basic fact:
\begin{equation} \label{basic1}
|| \eta|^s-|\eta'|^s| \,\leq \,|\eta-\eta'|^s \, , \quad 0<s<1
\end{equation}
for $\eta ,\eta'\in \R 2$, we obtain
\begin{displaymath}
|\hat{h}(\eta)| \,\leq\,\frac{1}{(2\pi)^2}\, 
\int_{\R 2} \,|\zeta|^s \,|\hat{f}(\zeta)|\ |\hat{g}(\eta-\zeta)|\,d\zeta \ .
\end{displaymath}
At last depending on how we apply Young's inequalities we obtain
either \eqref{commHs1} or \eqref{commHs2}.
\end{proof}

We now state the anounced estimates in Sobolev norms.

\begin{prop} Let $0 < s < 1$ and $\sigma >1$.\\
There exists $C > 0$ and $\kappa >0$ such that if $F$ belongs to
$H^s(\R2)$, $b$ to $L^{\infty} (\R2)$ with $\kappa \,
\norme{b}{\infty} < 1$ and $I^s\,b$ belongs to $H^{\sigma}(\R2)$, then, if $\Pi$ solves \eqref{pression-type}, $I^s\,\nabla \Pi$  belongs to
$L^2(\R2)$ and
\begin{equation} \label{pressionHs}
\norme{I^s\,\nabla \Pi}{2} 
\leq \frac{C}{1-\kappa \, \norme{b}{\infty}} \, 
( \norme{I^s\,F}{2} + \frac{1}{1-\kappa\,\norme{b}{\infty}} \, 
\hnorme{I^s\,b}{\sigma}\,\norme{F}{2}) \ .
\end{equation}
\end{prop}

\begin{proof}
Applying $I^s$ to \eqref{pression-type} and commuting $b$ and $I^s$ yields
\begin{displaymath}
\textnormal{div} \, \big((1+b) \nabla I^s\,\Pi \big) = 
\textnormal{div} \, \big([b,I^s] \nabla \Pi \big) + 
\textnormal{div} \, \big(I^s\,F \big) \, .
\end{displaymath}
Applying then \eqref{pressionLp} to this equation we obtain
\begin{displaymath}
\norme{I^s\,\nabla \Pi}{2} \,\leq\,
\frac{C}{1-\kappa \, \norme{b}{\infty}} \, 
(\norme{I^s\,F}{2}+\norme{[b,I^s]\nabla \Pi}{2}) \ .
\end{displaymath}
Now using first \eqref{commutateurHs2} then applying
\eqref{pressionLp} once again yields
\begin{displaymath}
\norme{[b,I^s]\nabla \Pi}{2} 
\,\leq\,C\,\hnorme{I^s\,b}{\sigma}\,\norme{\nabla \Pi}{2} 
\,\leq\,\frac{C}{1-\kappa\,\norme{b}{\infty}}\,
\hnorme{I^s\,b}{\sigma}\,\norme{F}{2} \ .
\end{displaymath}
Gathering everything leads to \eqref{pressionHs}.
\end{proof}

\section{Density equation}

In this section, we gather information on the following linearization
of the density equation:
\begin{equation} \label{lintransport}
\partial_t \, b + 
\big( \, ( \nu - \frac{1}{2} \, \xi ) \cdot \nabla \, \big) \, b =  0
\end{equation}
where $\widetilde{\nu}$ is a divergence-free vector-field, $\alpha \in
\R{}$ and $\nu = \alpha \, v^{G} + \widetilde{\nu}$. By linearization,
we mean that we do not assume that $\widetilde{\nu}$ is obtained from
a solution $\widetilde{w}$ of the vorticty equation in \eqref{complete},
 which involves $b$.

We begin with estimates in $L^p$-spaces, weighted or not. We recall that $L^p_w(\R2)$ for $1 \leq p \leq
\infty$ is the weighted space defined in \eqref{GLp}. The next proposition is the density part of the annouced
\emph{first step} of the proof.

\begin{prop}
Let $T > 0$. \\
Assume that $\widetilde{\nu}$ is a divergence-free vector field
belonging to $L^2 (0,T \, ; L^{\infty}(\R 2)) $. 
Then $b$, the solution of \eqref{lintransport} with initial data
$b_0$, satisfies
\begin{enumerate}
\item for any $1 \leq p \leq +\infty$, provided $b_0 \in L^p (\R 2)$, 
\begin{equation} \label{transportLp}
\norme{b(t)}{p}  \,  \leq \, \norme{b_0}{p} \, e^{-\frac tp} 
 \  , \quad \textrm{for} \  0 < t < T \  ;
\end{equation}

\item  for any $1 \leq p \leq +\infty$, provided $ G^{-\frac 12} \,
  b_0 \in L^p (\R{2})$,
\begin{equation} \label{transportGLp}
\gnorme{b( t )}{p}  \, \leq \, \gnorme{b_{0}}{p} \,  e^{-\frac{t}{p}} \, 
e^{\frac 18 \int_{0}^{t} \, \norme{\widetilde {\nu} (s) } {\infty} ^{2} \, ds}
 \  , \quad  \textrm{for} \  0 < t < T .
\end{equation}
\end{enumerate}
\end{prop}

\begin{proof}
\begin{enumerate}
\item For $1 \leq p < +\infty$, multiplying \eqref{lintransport} by $\textrm{sgn}(b)\,|b|^{p-1}$, where $\textrm{sgn}$ is the usual sign function, and integrating by parts yield
\begin{displaymath}
\frac{d}{dt}\,\norme{b}{p}^p  
=-\int_{\R2} \big(\,(\nu-\frac 12\,\xi)\cdot \nabla\,\big)\,|b|^p=-\norme{b}{p}^p
\end{displaymath}
since $\textnormal{div} \, \nu = 0$ and $\textnormal{div} \, \xi =
2$. Integrating gives \eqref{transportLp} in this case. The case
$p=+\infty$ follows letting $p$ go to infinity.

\item For $1 \leq p < +\infty$, starting as in the former case leads to
\begin{displaymath}
\frac{d}{dt} \, \gnorme{b}{p}^p
=-\int_{\R 2}  G^{-\frac p2} \big( \, (\nu-\frac 12 \, \xi) \cdot
\nabla \, \big) \, |b|^p
=\int_{\R 2}  G^{-\frac p2} \, \big( \, \frac p4 \, \xi \cdot (\nu-\frac 12
\, \xi)-1 \, \big) |b|^p
\end{displaymath}
since $\nabla G^{-\frac p2} =  G^{-\frac p2} \, \frac p4 \, \xi$. Now, since $\xi \perp v^{G}(\xi)$, we have
\begin{displaymath}
\xi \cdot(\nu(\xi)-\frac 12 \,\xi)
\,=\,\xi \cdot \widetilde{\nu}(\xi)-\frac 12 \,|\xi|^2
\,\leq\,\frac 12 \,|\widetilde{\nu}(\xi)|^2
\end{displaymath}
hence
\begin{displaymath}
\frac{d}{dt} \, \gnorme{b}{p}^p
\leq  (-1 + \frac p8 \, \norme{\widetilde{\nu}}{\infty}^2) \ 
\gnorme{b}{p}^p \  .
\end{displaymath}
Again integrating achieves the proof for finite $p$ and the case $p = + \infty$
follows letting $p$ go to infinity.
\end{enumerate}
\end{proof}

The next proposition corresponds to the density part of the
announced \emph{second step}: local-in-time estimates in Sobolev
norms. In order to prove some part of it we will need the following
commutator lemma.

\begin{lemme} Let $s \geq 1$ and $\sigma >1$.\\
There exists $C > 0$ such that if $I^{s} f$ belongs to 
$L^{2} (\R{2})$ and $g$ to $H^{\sigma}(\R{2})$, then $I^{s} (fg) - f
\, I^{s} g$ belongs to $L^{2} (\R{2})$ and 
\begin{equation}
\norme{I^{s} (fg) - f \, I^{s} g}{2} \,
\leq C \, \norme{I^{s} f}{2} \, \hnorme{g}{\sigma}
   + C \, \hnorme{\nabla f}{\sigma} \,  \norme{I^{s-1} g}{2}
\end{equation}
\end{lemme}

\begin{proof} The proof is essentially the same as the one of {\bf Lemma 1}
  except that here \eqref{basic1} is replaced by
\begin{equation} \label{basic2}
||\zeta|^s -|\zeta'|^s|\,\leq\,
C\,|\zeta-\zeta'|\,(|\zeta-\zeta'|^{s-1}+|\zeta'|^{s-1}), \quad \textrm{for} \ \zeta , \zeta' \in \R2 \ .
\end{equation}
Note that we could have obtained, as in {\bf Lemma 1}, various estimates depending on the way we apply Young's inequalities.
\end{proof}

\begin{prop}Let $T > 0$.
\begin{enumerate}
\item Let $ 0 < s \leq 2$ and $0 < \epsilon < s$.\\
Assume that $\widetilde{\nu}$ is a divergence-free vector field with
$\nabla \widetilde{\nu} \in L^1(0,T;H^1(\R2))$. \\ 
Then there exists $C_T>0$ independent of $\widetilde{\nu}$ such that, for any initial data $b_0 \in H^s(\R 2)$, any solution $b \in L^{\infty}(0,T;
H^{s-\epsilon}(\R2))$ of \eqref{lintransport} satisfies
\begin{equation} \label{transportHspetit}
\hnorme{b(t)}{s-\epsilon} \leq 
C_T\, \hnorme{b_0}{s}\,\exp \big((C_T \int_0^t\,
\hnorme{\nabla \nu(\tau)}{1}\,d\tau )^2 \big) \ , 
\quad \textrm{for} \ 0 \leq t \leq T.
\end{equation}

\item Let $s > 2$. \\
Assume that $\widetilde{\nu}$ is a divergence-free vector field whith
$\nabla \widetilde{\nu} \in L^1(0,T;H^{s-1}(\R 2))$. \\ 
Then \eqref{lintransport} has a unique solution $b \in L^{\infty}(0,T;
H^s(\R2))$, for any initial data $b_0 \in H^s(\R2)$. Moreover there 
exists $C > 0$ independent of $b_0$ and $\widetilde{\nu}$ such that 
$b$ satisfies
\begin{equation} \label{transportHsgrand}
\hnorme{b(t)}{s} \leq 
C\,\hnorme{b_0}{s}\,e^{\frac{s-1}{2} t}\,\exp \big(C \int_0^t\,\hnorme{\nabla  \nu(\tau)}{s-1}\,d\tau \big) \ ,\quad \textrm{for}\ 0 \leq t \leq T.
\end{equation}
\end{enumerate}
\end{prop}

\begin{proof}
\begin{enumerate}
\item See {\bf Theorem 0.1} in \cite{Danchin-transport}.

\item We can compute $[I^s,\frac{\xi}{2}]\cdot \,f=-\frac s2 \ I^{s-2} \ \textnormal{div}\,f$, for any vector field $f$. Thus applying $I^s$ to equation \eqref{lintransport} and commuting yield
\begin{displaymath}
\partial_t \,I^s b + 
\big(\,(\nu-\frac 12 \,\xi)\cdot \nabla\,\big)\,I^s b
=\frac s2 \,I^s b-[I^s,\nu] \cdot \,\nabla b \ .
\end{displaymath}
Then multiplying by $I^s b$ and integrating lead to
\begin{displaymath}
\frac 12 \,\frac{d}{dt}\,\norme{I^s b}{2}^2
-\frac{s-1}{2} \norme{I^s b}{2}^2 
=-\int_{\R2} I^s b\ [I^s,\nu]\cdot \,\nabla b
\end{displaymath}
since $\textnormal{div}\,\nu=0$. Now use Cauchy-Schwartz'
inequality and apply {\bf Lemma~2} (with $\sigma=s-1$) to get
\begin{displaymath}
\frac 12 \,\frac{d}{dt}\,\norme{I^s b}{2}^2
-\frac{s-1}{2}\,\norme{I^s b}{2}^2\,
\leq\,C\,\hnorme{\nabla \nu}{s-1}\,\hnorme{b}{s}^2\ .
\end{displaymath}
At last combine the former with $\frac 12 \,\frac{d}{dt}\,\norme{b}{2}^2+\frac 12 \,\norme{b}{2}^2 \leq 0$ to obtain
\begin{displaymath}
\frac 12 \,\frac{d}{dt}\,\hnorme{b}{s}^2
-\frac{s-1}{2}\,\hnorme{b}{s}^2 
\,\leq\, C\,\hnorme{\nabla \nu}{s-1}\,\hnorme{b}{s}^2
\end{displaymath}
which yields \eqref{transportHsgrand} by a mere integration.
\end{enumerate}
\end{proof}

The last estimate we state for the linearized transport equation \eqref{lintransport} is intended to be used for the proofs of the convergence of our iterative scheme and of the uniqueness of our solutions. Indeed we estimate the difference of two solutions of equations of type \eqref{lintransport}.

\begin{prop}
Let  $T>0$. \\
Assume that, for $i=1,2$, $\widetilde{\nu_{i}}$ is a divergence-free
vector field belonging to $L^{2} (0,T; W^{1,\infty}(\R{2})) $.
If, for $i=1,2$, $b_{i}$ is a solution of 
\begin{displaymath}
\partial_t \, b_i + 
\big( \, (\nu_i-\frac 12 \, \xi) \cdot \nabla \, \big) \, b_i= 0  \  ,
\end{displaymath}
where $\nu_i = \alpha \, v^G + \widetilde{\nu}_i$, 
with initial data $b_0$, then $b_1$ and $b_2$ satisfy 
\begin{enumerate}
\item provided that $G^{-\frac 12} \nabla b_0$ belongs to 
$L^p (\R 2)$, for some $1 \leq p \leq + \infty$, 
\begin{equation} \label{lemmetransportdelta} 
\gnorme{\nabla b_i(t)}{p} \leq \gnorme{\nabla b_0}{p} \, 
e^{-t (\frac 1p - \frac 12)}  \, 
e^{\frac 18 \int_0^t \, \norme{\widetilde{\nu}_i(s)}{\infty}^2 \, ds} 
\, e^{\int_0^t \, \norme{\nabla \nu_i(s)}{\infty} \, ds}
\end{equation}
for $i=1,2$ and $0 \leq t \leq T$;

\item provided that $b_0$ is such that $G^{-\frac 12} b_0$ belongs to $L^p (\R 2)$ and $G^{-\frac 12} \nabla b_0$ to $L^q (\R 2)$, for some $1 \leq p < q \leq +\infty$, 
\begin{equation} \label{transportdelta}
\gnorme{(b_2-b_1)(t)}{p} \leq 
e^{\frac 18 \int_0^t \, \norme{\widetilde{\nu}_2(s)}{\infty}^2 \, ds}
\, \sup_{0 \leq s \leq t} \gnorme{\nabla b_1(s)}{q}
\, \int_0^t \, \norme{(\widetilde{\nu}_2 - \widetilde{\nu}_1)(s)}{r} \, ds
\end{equation}
for $0 \leq t \leq T$, where $r$ is such that 
$\frac 1p = \frac 1q + \frac 1r$.
\end{enumerate}
\end{prop}

\begin{proof}
\begin{enumerate}
\item Derive the equation for $b_1$ to get for $j=1,2$,
\begin{displaymath}
\partial_t \, \partial_j b_1 + 
\big( \, (\nu_1-\frac 12 \, \xi) \cdot \nabla \, \big) \, \partial_j b_1 
=-\partial_j \nu_1 \cdot \nabla b_1 + \frac 12 \, \partial_j b_1 \ .
\end{displaymath}
From this, following the proof of \eqref{transportGLp}, we obtain for $j=1,2$, 
\begin{displaymath}
\frac{d}{dt} \gnorme{\partial_j b_1(t)}{p}^p \leq 
(-1 + \frac p2 +\frac p8 \, \norme{\widetilde{\nu}_1(t)}{\infty}^2)
\, \gnorme{\partial_j b_1(t)}{p}^p 
+ p \, \norme{\nabla \nu_1(t)}{\infty} \, \gnorme{\nabla b_1(t)}{p}^p
\ .
\end{displaymath}
Now combining the inequalities for $j=1,2$ and integrating lead to \eqref{lemmetransportdelta}.

\item Observe that $b_2-b_1$ satisfies 
\begin{displaymath}
\partial_t \, (b_2 - b_1) + 
\big( \, (\nu_2 - \frac 12 \, \xi) \cdot \nabla \, \big) \, (b_2 - b_1) 
= - (\widetilde{\nu}_2 - \widetilde{\nu}_1) \cdot \nabla b_1 \ .
\end{displaymath}
Now following again the proof of \eqref{transportGLp} yields
\begin{displaymath}
\frac{d}{dt} \, \gnorme{\delta \, b}{p}^p \,
\leq \, (-1 + \frac p8 \, \norme{\widetilde{\nu}_2}{\infty}^2) \
\gnorme{\delta \, b}{p}^p + p \, \gnorme{\delta \, b}{p}^{p-1} \,
\gnorme{\nabla  b_1}{q} \, \norme{\widetilde{\nu}_2-\widetilde{\nu}_1}{r}
\end{displaymath}
where $\delta \, b = b_2 -b_1$. It is now straightforward to derive \eqref{transportdelta}.
\end{enumerate}
\end{proof}

\section{Vorticity equation}

As announced, we now study a linearization of the vorticty equation:
\begin{equation} \label{vorticitelin}
\partial_t \, \widetilde{w} 
      - ( \mathcal{L} - \alpha \, \Lambda ) \, \widetilde{w} 
      + \big(\, \widetilde{\nu} \cdot \nabla \, \big) \, \widetilde{w} 
=  \textnormal{div} \, \big( b \, (\nabla w
                                          + \nabla^{\perp} \Pi ) \big)
\end{equation}
where $\mathcal{L}$ and $\Lambda$ are as in \eqref{LandLambda}, $b$ is
a real function, $\widetilde{\nu}$ is a divergence-free vector field,
$\alpha \in \R{}$,
\begin{displaymath}
\begin{array}{rclcrcl}
\widetilde{v} & = & K_{BS} \star \widetilde{w} & , &
v & = & \alpha \, v^G + \widetilde{v} \ ,\\
w & = & \alpha \, G + \widetilde{w} & , &
\nu & = & \alpha \, v^G + \widetilde{\nu} \ ,
\end{array}
\end{displaymath}
and $\nabla \Pi$ is obtained by solving
\begin{equation} \label{pressionlin}
\textnormal{div} \, \big((1+b) \, \nabla \Pi \big)  
=\textnormal{div} \, \big((1 +b) \, \triangle v 
                         - (\nu \cdot \nabla) \, v \big) \ .
\end{equation}

Remind that we always assume $\int_{\R 2} \widetilde{w}_0 = 0$.

\subsection{Global estimate}

In this subection we establish a global-in-time estimate in weighted
$L^p$-spaces.

\begin{prop} \label{premieresestimations}
Let $\alpha \in \R{}$, $K_0 > 0$. There exist $\varepsilon_0 > 0$ and 
$C >0$ such that if $b$ is a real function and $\widetilde{\nu}$ a
divergence-free vector field such that
\begin{enumerate} 
\item for $0 < t < T$, for any $1 \leq p \leq +\infty$, $2\leq q \leq +\infty$,
\begin{displaymath}
\begin{array}{lcrclcr}
\norme{b(t)}{p}  &\leq & \norme{b_0}{p} \, e^{-\frac tp} &,&
\gnorme{b(t)}{q} &\leq & \gnorme{b_0}{q} \,  e^{-\frac tq} \, e^{K_0} 
\end{array}
\end{displaymath}

\item for $0 < t < T$,
\begin{displaymath}
\begin{array}{lcrclcr}
\norme{\widetilde{\nu}(t)}{8} & \leq & K_0 &,&
\int_0^t \norme{\widetilde{\nu}}{\infty}^2 & \leq & \frac{1}{24} 
\end{array}
\end{displaymath}

\item and
\begin{displaymath}
\begin{array}{lcrclcr}
\gnorme{b_0}{4} & \leq & \varepsilon_0 &,&
\gnorme{b_0}{\infty} & \leq & \varepsilon_0
\end{array}
\end{displaymath}
\end{enumerate}
then any solution $\widetilde{w} \in L^{\infty} (0,T \, ; L^2_w (\R 2))$ of \eqref{vorticitelin}, with initial data $\widetilde{w}_0 \in L^2_w (\R 2)$, satisfies, for any $0<t<T$,
\begin{eqnarray} \label{vorticitelinGLp}
\gnorme{\widetilde{w}(t)}{2} ^2
 &+& C \int_0^t \, \big( \gnorme{\widetilde{w}}{2}^2 
 + \gnorme{\nabla \widetilde{w}}{2}^2
 + \gnorme{|\xi| \widetilde{w}}{2}^2 \big)\\
 &\leq& 2 \, \gnorme{\widetilde{w}_0}{2}^2 \ 
 + C \, |\alpha | \, \gnorme{b_0}{4}
\, .
\end{eqnarray}
\end{prop}

Note that the assumptions on $b$ corresponds to the first proposition
of the previous section. Note also that once $\alpha$ and $K_0$ are fixed,
since $L^2_w$ is embedded in any $L^p$, $1 \leq p \leq 2$, and $H^1$
is embedded in any $L^q$, $2 \leq q<\infty$, inequalities
\eqref{BiotSavartL8} and \eqref{vorticitelinGLp} enable us to make $\int_0^t
\norme{\widetilde{v}}{\infty} ^2$ as small as we want provided we take
$\widetilde{w}_0$ and $b_0$ small enough. At last note that the proposition enables us to bound $\int_0^t \hnorme{\nabla \widetilde{v}}{1}^2$, which will be used in \eqref{transportHspetit}.

\begin{proof}
Our strategy is to multiply \eqref{vorticitelin} by $G^{-1} \,
\widetilde{w}$ and integrate to bound $\frac{d}{dt}
\gnorme{\widetilde{w}}{2}^2$. In what follows, we examine each term arising once multiplyed by $G^{-1}\,\widetilde{w}$ and integrated on $\R2$.

\begin{itemize}
\item   Let us emphasize first that \eqref{vorticitelin} preserves
$\int_{\R 2} \, \widetilde{w}$. Hence $\int_{\R 2} \, \widetilde{w}=0$.

Let $\textnormal{L}:=G^{-\frac 12}\,(-\mathcal{L})\,G^{\frac 12}$. A direct calculation shows that $\textnormal{L}=-\triangle+~\frac{|\xi|^2}{16}-~\frac12$ is a harmonic oscillator with spectrum $\{0,\frac12,1,\frac32,\ldots\}$. Moreover $0$ is a simple eigenvalue with eigenvector $G^{\frac12}$. In particular, if $f$ belongs to the domain of $\textnormal{L}$ with $\int_{\R2} G^{\frac 12}\,f=0$, then $\int_{\R2} f\,\textnormal{L}f \geq \frac 12 \,\norme{f}{2}^2$. 

Coming back to $\mathcal{L}$, we obtain: if $G ^{-\frac 12}
\widetilde{w}$ belongs to the domain of $\textnormal{L}$ with
$\int_{\R 2} \, \widetilde{w}=0$, then, for any $0 < \gamma < \frac 12$, 
\begin{displaymath}
\int_{\R 2} G^{-1} \, \widetilde{w} \  \mathcal{L} \, \widetilde{w}
\, \leq \,  -\frac 12 \, (1 - \gamma) \, \gnorme{\widetilde{w}}{2}^2
+ \gamma \, \int_{\R 2} G^{-1}\,\widetilde{w} \, \mathcal{L} \widetilde{w}
\end{displaymath}
thus integrating by part from the formula for $\textnormal{L}$
\begin{displaymath}
\int_{\R 2} G^{-1} \, \widetilde{w} \  \mathcal{L} \, \widetilde{w}
\, \leq \,  -\frac 12 \, (1 - 2 \gamma) \, \gnorme{\widetilde{w}}{2}^2 
- \gamma \, \big( \norme{\nabla (G^{-\frac 12} \, \widetilde{w})}{2}^2
+ \gnorme{\frac{|\xi|}{4} \, \widetilde{w}}{2} ^{2} \big)
\end{displaymath}
and expanding
\begin{displaymath} \begin{array}{rcl}
\int_{\R 2} G^{-1} \, \widetilde{w} \  \mathcal{L} \, \widetilde{w}   
&\leq &-\frac 12 \  (1 - 2 \gamma) \  \gnorme{\widetilde{w}}{2}^2 \\[1ex]
& & -\gamma \  \big( \, \gnorme{\nabla \widetilde{w}}{2}^2 
+ 2 \, \gnorme{\frac{|\xi|}{4} \, \widetilde{w}}{2}^2 
+ 2 \int_{\R 2} G^{-1} \, \nabla \widetilde{w} \cdot \frac{\xi}{4} \,
\widetilde{w} \, \big) 
\end{array}
\end{displaymath}
hence 
\begin{equation} \label{oscillateur}
\int_{\R 2} G^{-1} \, \widetilde{w} \  \mathcal{L} \, \widetilde{w} 
 \leq -\frac 12 \, (1 - 2 \gamma) \, \gnorme{\widetilde{w}}{2}^2
- \gamma \, \big(\frac 13 \, \gnorme{\nabla \widetilde{w}}{2}^2 
+ \frac 12 \, \gnorme{\frac{|\xi|}{4} \, \widetilde{w}}{2}^2 \big) \, .
\end{equation}

\item Recalling that $\Lambda \widetilde{w}=v^G \cdot \nabla \widetilde{w} + \widetilde{v} \cdot \nabla G$, we obtain
\begin{equation}
\int_{\R 2} G^{-1} \, \widetilde{w} \, \Lambda \widetilde{w} = 0 \ .
\end{equation}
Indeed , from $v^G(\xi) \perp \xi$ and $\nabla G^{-1}=-\frac{\xi}{2}G^{-1}$, we derive
\begin{displaymath}
\int_{\R 2} G^{-1} \,\widetilde{w}\,v^{G}\cdot \nabla \widetilde{w}=-\frac 12 \int_{\R 2} G^{-1}\,\frac{\xi}{2}\,\cdot\,v^{G}\,\widetilde{w}^2 = 0
\end{displaymath}
And, on the other hand, using the identity $\eta^{\perp} \cdot \xi = - \xi^{\perp} \cdot \eta$ and the explicit formula \eqref{BiotSavart} for the
Biot-Savart law, we derive 
\begin{displaymath}
\begin{array}{rcl}
\int_{\R 2} G^{-1}\,\widetilde{w}\ \widetilde{v}\cdot \nabla G
 &=&-\int_{\R 2} \widetilde{w}(\xi)\ \widetilde{v}(\xi)\cdot
\frac{\xi}{2}\,d\xi \\[1ex]
 &=&-\frac{1}{4 \pi} \IInt{\R 2 \times \R 2} 
\widetilde{w}(\xi)\ \frac{(\xi-\eta)^{\perp}\cdot\,\xi}{|\xi-\eta|^2}
\ \widetilde{w}(\eta)\,d\eta\,d\xi \\[1ex]
 &=&\frac{1}{4 \pi} \IInt{R 2 \times \R 2}
\widetilde{w}(\xi)\ \frac{\eta^{\perp}\cdot\,\xi}{|\xi-\eta|^2}
\ \widetilde{w}(\eta)\,d\eta\,d\xi \\[1ex]
 &=&-\frac{1}{4 \pi} \IInt{\R 2 \times \R 2}
\widetilde{w} (\xi)\ \frac{\xi^{\perp} \cdot \eta}{|\xi-\eta|^2}
\ \widetilde{w}(\eta)\,d\eta\,d\xi \\[1ex]
 &=&-\int_{\R 2} G^{-1}\,\widetilde{w}\ \widetilde{v}\cdot \nabla G  
\ .
\end{array}
\end{displaymath}
Thus $\int_{\R 2} G^{-1}\,\widetilde{w}\ \widetilde{v}\cdot\nabla G = 0$.

\item Using H\"older's inequalities, we also obtain
\begin{equation}
|\int_{\R 2} G^{-1} \, \widetilde{w} \  \widetilde{\nu} \cdot \nabla
\widetilde{w}| 
\, \leq \, 6 \, \norme{\widetilde{\nu}}{\infty} ^{2} \,
\gnorme{\widetilde{w}}{2}^2  + \frac{1}{24} \, \gnorme{\nabla
  \widetilde{w}}{2}^2 \ .
\end{equation}

\item Integrating by part, we obtain
\begin{displaymath}
\int_{\R 2} G^{-1} \, \widetilde{w} \  \textnormal{div} \, (b \, \nabla
\widetilde{w})
=-\int_{\R 2}  G^{-1} \, b \, |\nabla \widetilde{w}|^2 
- \frac 12 \int_{\R 2}  G^{-1} \, b \, \widetilde{w} \, \xi \cdot
\nabla \widetilde{w} 
\end{displaymath}
then, using H\"older's inequalities and the fact that $\norme{b(t)}{\infty} \leq \norme{b_0}{\infty}$,
\begin{equation}
|\int_{\R2} G^{-1}\,\widetilde{w}\ \textnormal{div}\,(b\,\nabla \widetilde{w})|\,
\leq \, \frac 54 \, \norme{b_0}{\infty} \, \gnorme{\nabla \widetilde{w}}{2}^2 
+ \frac 14 \, \norme{b_0}{\infty} \, \gnorme{|\xi| \, \widetilde{w}}{2}^2 
\ .
\end{equation}

\item
In the same way, since $\norme{b(t)}{2} \leq \norme{b_0}{2}e^{-\frac t2}$, we have
\begin{displaymath} \begin{array}{rcl}
|\int_{\R 2} G^{-1} \, \widetilde{w} \ \textnormal{div} \, (b \, \nabla G)| 
&=& |\frac 12 \int_{\R 2} b \, \xi \, \cdot \, \nabla \widetilde{w} 
+ \frac 14 \, \int_{\R 2} b \, \widetilde{w} \, |\xi|^2 | \\[1ex]
&\leq& C \, \norme{b_0}{2} \, e^{- \frac t2} \,
( \gnorme{\nabla \widetilde{w}}{2} + \gnorme{\widetilde{w}}{2} )
\end{array}
\end{displaymath}
thus
\begin{equation}
|\int_{\R 2} G^{-1} \, \widetilde{w} \ \textnormal{div} \, (b \, \nabla G)| \,
\leq \, C \, \norme{b_0}{2} \, (e^{-t} + \gnorme{\nabla
  \widetilde{w}}{2}^2 + \gnorme{\widetilde{w}}{2}^2)  
\ .
\end{equation}

\item Finally, integrating by part, using H\"older's inequalities and
  applying inequality \eqref{pressionLp}, we obtain, for $b_0$ small
  enough in $L^{\infty}$, 
\begin{eqnarray*} 
|\int_{\R 2} G^{-1}\,\widetilde{w}\ \textnormal{div}\,(b\,\nabla^{\perp} \Pi)| 
&=&|\frac 12\,\int G^{-1}\,b\,\widetilde{w}\,\xi \cdot \nabla^{\perp}\Pi 
+\int_{\R 2} G^{-1}\,b\,\nabla \widetilde{w} \cdot \nabla^{\perp}\Pi| \\
&\leq& \frac{C}{1-\kappa\,\norme{b_0}{\infty}}\, 
(\gnorme{|\xi| \widetilde{w}}{2}+\gnorme{\nabla \widetilde{w}}{2})\\
&\times&(\gnorme{b_0}{\infty}\,e^{K_0}\, 
\norme{(1+b)\,\triangle \widetilde{v}}{2}\\
& &+\,\gnorme{b_0}{4}\,e^{K_0}\,e^{-\frac t4} 
\norme{\alpha \,(1+b)\,\triangle v^G-(\nu \cdot \nabla)\,v}{4}) \ .
\end{eqnarray*}
Now, on one hand,since $\norme{b(t)}{\infty} \leq \norme{b_0}{\infty}$, estimate \eqref{BiotSavartHsGrand} yields
\begin{displaymath}
\norme{(1+b)\,\triangle \widetilde{v}}{2} \leq C\,(1+\norme{b_0}{\infty})\,
\norme{\nabla \widetilde{w}}{2}\ .
\end{displaymath}
On the other hand, similarly, we have
\begin{displaymath}
\norme{(1+b)\,\triangle v^G}{4} \leq C\,(1+\norme{b_0}{\infty})\ .
\end{displaymath}
At last, H\"older's inequalities, estimate \eqref{BiotSavartGrad} and Sobolev's embeddings yields
\begin{displaymath}
\norme{(\nu \cdot \nabla)\,v}{4} \leq C\,(|\alpha|+\norme{\widetilde{\nu}}{8})\,(|\alpha|+\norme{\widetilde{w}}{8}) \leq C\,(|\alpha|+\norme{\widetilde{\nu}}{8})\,(|\alpha|+\hnorme{\widetilde{w}}{1}) \ .
\end{displaymath}
Taking this into account yields when $\kappa\,\norme{b_0}{\infty} \leq \frac 12$
\begin{eqnarray}
|\int_{\R 2} G^{-1} \, \widetilde{w}
&\times& \textnormal{div}\,(b\,\nabla^{\perp} \Pi)|
\nonumber \\
 &\leq& e^{- \frac t2}\qquad \quad
C\,|\alpha|\,e^{K_0}\,\gnorme{b_0}{4}\,(1+|\alpha|+\norme{\widetilde{\nu}}{8})^2
\nonumber \\
 &+& \gnorme{\widetilde{w}}{2}^2\quad \ \ 
C\,e^{K_0}\,\gnorme{b_0}{4}\,(1+|\alpha|+\norme{\widetilde{\nu}}{8})^2 
\nonumber \\
 &+& \gnorme{\nabla \widetilde{w}}{2}^2 \ \ \,\,C\,e^{K_0}(\gnorme{b_0}{\infty}
+\gnorme{b_0}{4}(1 +|\alpha|+\norme{\widetilde{\nu}}{8}))
\nonumber \\
 &+& \gnorme{|\xi| \widetilde{w}}{2}^2\ \ C\,e^{K_0}\,(\gnorme{b_0}{\infty}
+\gnorme{b_0}{4}(1+|\alpha|)) 
\ .
\end{eqnarray}
\end{itemize}

It only remains to us to gather everything after setting $\gamma=\frac 14$ in \eqref{oscillateur} and integrate in time in order to obtain, when $\kappa\,\norme{b_0}{\infty} \leq \frac 12$,
\begin{displaymath} \begin{array}{lllll}
\gnorme{\widetilde{w}(t)}{2}^2 &\times&
(1-12\int_0^t \norme{\widetilde{\nu}}{\infty}^2)& & \\[1ex]
 &+&C\,\int_0^t \,\gnorme{\widetilde{w}}{2}^2 &\times&
(1-e^{K_0} \gnorme{b_0}{4} (1+|\alpha|+\norme{\widetilde{\nu}}{8})^2)\\[1ex]
 &+&C\,\int_0^t \,\gnorme{\nabla \widetilde{w}}{2}^2 &\times&
(1-\,e^{K_0}(\gnorme{b_0}{\infty}+\gnorme{b_0}{4}\,
(1+|\alpha|+\norme{\widetilde{\nu}}{8})) \\[1ex]
 &+&C\,\int_0^t \gnorme{|\xi|\,\widetilde{w}}{2}^2 &\times&
(1-e^{K_0}(\gnorme{b_0}{\infty}+\gnorme{b_0}{4}\,(1+|\alpha|))) \\[1ex]
 &\leq& \gnorme{\widetilde{w}_0}{2}^2\ 
+\ C\,|\alpha | &\times& \gnorme{b_0}{4}\,e^{K_0}\,
(1+|\alpha|+\norme{\widetilde{\nu}}{8})^2)
\end{array}
\end{displaymath}
which yields the proposition since $\int_0^t
\norme{\widetilde{\nu}}{\infty}^2 \, \leq \, \frac{1}{24}$.
\end{proof}

\subsection{Sobolev estimate}

In this subsection we prove a local-in-time estimate in Sobolev norms
for solutions of equation \eqref{vorticitelin}. Remind that 
$I=(-\triangle)^{\frac 12}$.

\begin{prop} $ $ \\
Let $0 < s < 1$, $1+s < s' < 2$, $1 < s''< 2-s$ and $\alpha \in \R{}$.\\
There exists $\varepsilon_0>0$ and, for $K>0$, there exists $C >0$ such that if $b$ is a real function and $\widetilde{\nu}$ a divergence-free vector field such that
\begin{displaymath}
\begin{array}{lllclllc}
\sup_{[0,T]}\norme{b}{\infty}&\leq&\norme{b_0}{\infty} \leq \varepsilon_0
&,& \sup_{[0,T]}\hnorme{b}{s'}&\leq&K &, \\[1ex]
\sup_{[0,T]}\norme{\widetilde{w}}{2}&\leq&K
&,& \int_0^T \gnorme{\widetilde{w}}{2}^2&\leq&K &, \\[1ex]
\int_0^T \norme{\nabla \widetilde{w}}{2}^2&\leq&K
&,& \int_0^T \norme{\widetilde{\nu}}{\infty}^2&\leq&K &, \\[1ex]
\sup_{[0,T]} \norme{I^s \widetilde{\nu}}{2}&\leq&\varepsilon_0
&,& \int_0^T \hnorme{I^s \widetilde{\nu}}{s''}^2&\leq&K &,
\end{array}
\end{displaymath}
then any solution $\widetilde{w} \in L^{\infty}(0,T;\dot{H}^s(\R 2))$ of 
\eqref{vorticitelin}, with initial data $\widetilde{w}_0 \in \dot{H}^s (\R 2)$,
 satisfies for any $0 < t < T$,
\begin{equation}\label{vorticitelinHs}
\norme{I^s \widetilde{w} \, (t)}{2}^2 
+ C \, \int_0^t \norme{I^s \nabla \widetilde{w}}{2}^2 \ 
\leq \ C \ e^{C t} \ (\, \norme{I^s \widetilde{w}_0}{2}^2 + K \,)\ .
\end{equation}
\end{prop}

Note that \eqref{transportLp} and \eqref{transportHspetit} could provide us the validity of the asumptions on $b$, and \eqref{vorticitelinGLp} both the validity of the asumptions on $\widetilde{w}$ and, thanks to estimates \eqref{BiotSavartL8}, \eqref{BiotSavartHsGrand} and \eqref{BiotSavartHsPetit}, the validity of estimates on $\widetilde{\nu}$ when $\widetilde{\nu}=\widetilde{v}$. Conversely \eqref{vorticitelinHs} could be used in \eqref{transportHsgrand}, again when $\widetilde{\nu}=\widetilde{v}$.

\begin{proof}
We choose $\sigma$ such that $1 < \sigma < 1+s$, $s+\sigma < s'$ and $\sigma < s''$. The only role of $\sigma$ is to make clearer our use of commutator estimates \eqref{commutateurHs1} and \eqref{commutateurHs2}.

Our strategy is to apply $I^s$ to \eqref{vorticitelin}, then multiply by $I^s \widetilde{w}$ and estimate each term arising, in order to bound $\frac{d}{dt} \hnorme{I^s \widetilde{w}}{2}^2$.

\begin{itemize}
\item First we compute the commutator $ [I^s,\mathcal{L}] = \frac s2 I^s $ and obtain
\begin{equation}
\int_{\R 2} I^s \widetilde{w}\quad I^s \mathcal{L} \widetilde{w} =
-\int_{\R 2} |\nabla I^s \widetilde{w}|^2 
+ \frac{1+s}{2} \int_{\R 2} |I^s \widetilde{w}|^2 
\end{equation}
since integrating by parts yields 
$\int_{\R2} f\ \mathcal{L}f=-\int_{\R2}|\nabla f|^2+\frac 12 \int_{\R2}|f|^2$.

\item On one hand, we have
\begin{displaymath}
\int_{\R 2} I^s \widetilde{w} \quad I^s ((v^G \cdot \nabla) \, \widetilde{w})
=\int_{\R 2} I^s \widetilde{w} \quad (v^G \cdot \nabla) \, I^s \widetilde{w} \,
+\int_{\R 2} I^s \widetilde{w} \quad ([I^s,v^G] \cdot \nabla)\, \widetilde{w}
\end{displaymath}
with, integrating by parts,
\begin{displaymath}
\int_{\R 2} I^s \widetilde{w} \quad (v^G \cdot \nabla)\,I^s \widetilde{w}\,
=\,\frac 12 \int_{\R 2} (v^G \cdot \nabla)\,|I^s \widetilde{w}|^2\,=\,0
\end{displaymath}
and, using H\"older's inequalities and inequality \eqref{commutateurHs2},
\begin{displaymath}
|\int_{\R 2} I^s \widetilde{w} \quad ([I^s,v^G] \cdot \nabla)\,\widetilde{w}|\ \leq\ C\,\norme{I^s \widetilde{w}}{2}\,\hnorme{I^s v^G}{\sigma}\,
\norme{\nabla \widetilde{w}}{2} \ .
\end{displaymath}

On the other hand, we have
\begin{displaymath}
\int_{\R 2} I^s \widetilde{w} \quad I^s ((\widetilde{v} \cdot \nabla)\,G) 
=\int_{\R 2} I^s \widetilde{w} \quad (\widetilde{v} \cdot \nabla)\,I^s G
+\int_{\R 2} I^s \widetilde{w} \quad ([I^s,\widetilde{v}] \cdot \nabla)\,G
\end{displaymath}
with, using H\"older's inequalities,
\begin{displaymath}
|\int_{\R 2} I^s \widetilde{w} \quad (\widetilde{v} \cdot \nabla)\,I^s G|
\leq C\,\norme{I^s \widetilde{w}}{2} \norme{\widetilde{v}}{\infty} 
\norme{\nabla I^s G}{2}
\end{displaymath}
and, using both H\"older's inequalities and inequality \eqref{commutateurHs1},
\begin{displaymath}
|\int_{\R 2} I^s \widetilde{w} \quad ([I^s, \widetilde{v}] \cdot \nabla)\,G|
\leq C\,\norme{I^s \widetilde{w}}{2}\,\norme{I^s \widetilde{v}}{2}\,
\hnorme{\nabla G}{\sigma}\, .
\end{displaymath}

Using estimate \eqref{BiotSavartHsPetit} and estimate \eqref{BiotSavartL8} combined with Sobolev's embeddings, we derive
\begin{equation}
|\int_{\R2} I^s \widetilde{w} \quad I^s \Lambda \widetilde{w}|
\,\leq\,C\,(\gnorme{\widetilde{w}}{2}^2+\norme{\nabla \widetilde{w}}{2}^2)\ .
\end{equation}

\item In the same way, we have
\begin{displaymath}
\int_{\R2} I^s \widetilde{w} \quad (\widetilde{\nu} \cdot \nabla)\,I^s
\widetilde{w}=0
\end{displaymath}
and, applying H\"older's inequalities and inequality \eqref{commutateurHs1},
\begin{displaymath}
|\int_{\R2} I^s \nabla \widetilde{w} \cdot([I^s,\widetilde{\nu}]\,\widetilde{w})|
\,\leq\,C\,\norme{I^s \nabla \widetilde{w}}{2}\,\norme{I^s \widetilde{\nu}}{2}\,\hnorme{\widetilde{w}}{\sigma}
\end{displaymath}
thus, since $0\leq \sigma \leq 1+s$,
\begin{equation}
|\int_{\R2} I^s \widetilde{w} \quad I^s((\widetilde{\nu} \cdot \nabla)\,\widetilde{w})| \,\leq\,
C\,\norme{I^s \widetilde{\nu}}{2}\,\norme{I^s \nabla \widetilde{w}}{2}^2\,
+\,C\,\norme{I^s \widetilde{\nu}}{2}\,\norme{\widetilde{w}}{2}^2\ .
\end{equation}

\item Integrating by parts yields
\begin{displaymath}
\int_{\R2} I^s \widetilde{w} \quad I^s \textnormal{div}(\,b\nabla \widetilde{w})=-\int_{\R2} b\,|I^s \nabla \widetilde{w}|^2
-\int_{\R2} I^s \nabla \widetilde{w} \quad [I^s,b]\,\nabla \widetilde{w}\ .
\end{displaymath}
Applying \eqref{commutateurHs2}, we derive
\begin{equation}
|\int_{\R2} I^s\widetilde{w}\quad I^s \textnormal{div}(\,b\nabla \widetilde{w})| \,\leq\,(\norme{b_0}{\infty}+\varepsilon)\,\norme{I^s\nabla \widetilde{w}}{2}^2 
+\frac{C}{\varepsilon}\,\hnorme{b}{s+\sigma}^2\,\norme{\nabla \widetilde{w}}{2}^2
\end{equation}
where $\varepsilon>0$ is intended to be chosen small enough.

\item Similarly,
\begin{equation}
|\int_{\R2} I^s \widetilde{w} \quad I^s \textnormal{div}(\,b\nabla G)|\, 
\leq\,(\norme{b_0}{\infty}+\varepsilon)\,\norme{I^s \nabla \widetilde{w}}{2}^2+ \,C\,(\norme{b_0}{\infty}+\frac{1}{\varepsilon}\,\hnorme{b}{s+\sigma}^2)\ .
\end{equation}

\item First, integrating by parts yields
\begin{displaymath}
\int_{\R2} I^s \widetilde{w} \quad I^s \textrm{div}(b\nabla^{\perp} \Pi)
=-\int_{\R2} I^s \nabla \widetilde{w}\ \cdot\ b \,I^s \nabla^{\perp} \Pi
-\int_{\R 2} I^s \nabla \widetilde{w}\ \cdot\ ([I^s,b] \nabla^{\perp} \Pi)
\end{displaymath}
with, using H\"older's inequalities,
\begin{displaymath}
|\int_{\R2} I^s \nabla \widetilde{w} \ \cdot\ b\,I^s \nabla^{\perp} \Pi| 
\,\leq\,\norme{I^s \nabla \widetilde{w}}{2}\,\norme{b_0}{\infty} 
\,\norme{I^s \nabla \Pi}{2}
\end{displaymath}
and, using H\"older's inequalities and inequality \eqref{commutateurHs2},
\begin{displaymath}
|\int_{\R 2} I^s \nabla \widetilde{w} \cdot ([I^s,b] \nabla^{\perp} \Pi)|
\, \leq \, C \, \norme{I^s \nabla \widetilde{w}}{2} \, \hnorme{b}{s+\sigma} 
\, \norme{\nabla \Pi}{2} \ .
\end{displaymath}

Besides, on one hand, estimate \eqref{pressionLp} applied to equation \eqref{pressionlin} and H\"older's inequalities imply
\begin{displaymath}
\norme{\nabla \Pi}{2} \, \leq \, C \, (|\alpha| + |\alpha|^2 + 
\norme{\triangle \widetilde{v}}{2} + |\alpha| \norme{\nabla \widetilde{v}}{2}
+ \norme{\widetilde{\nu}}{\infty} \norme{\nabla \widetilde{v}}{2}
+ |\alpha| \norme{\widetilde{\nu}}{\infty})
\end{displaymath}
thus with estimates on Biot-Savart kernel
\begin{displaymath}
\norme{\nabla \Pi}{2} \, \leq \, C \, (|\alpha| + |\alpha| \norme{\widetilde{\nu}}{\infty} +|\alpha|^2 + (|\alpha| +\norme{\widetilde{\nu}}{\infty}) \, \norme{\widetilde{w}}{2} + \norme{\nabla \widetilde{w}}{2}) \ .
\end{displaymath}
On the other hand, estimate \eqref{pressionHs} applied to \eqref{pressionlin} and H\"older's inequalities imply
\begin{displaymath}
\norme{I^s \nabla \Pi}{2} \, \leq \, 
C \, (\hnorme{b}{s+\sigma} \norme{\nabla \Pi}{2}
+ \norme{I^s ((1+b) \, \triangle v)}{2}
+ \norme{I^s (\nu \cdot \nabla) \, v}{2})
\end{displaymath}
with, commuting $I^s$ and $b$ thanks to \eqref{commutateurHs2}, after some calculation,
\begin{displaymath}
\norme{I^s ((1+b) \, \triangle v)}{2}
\, \leq \, C \, \hnorme{b}{s+\sigma} \, 
(|\alpha| + \norme{\nabla \widetilde{w}}{2})
+ C \, (1+\norme{b_0}{\infty}) \, 
(|\alpha| + \norme{I^s \nabla \widetilde{w}}{2})
\end{displaymath}
and, in the same way, commuting $I^s$ and $\nu$,
\begin{displaymath}
\norme{I^s ((\nu \cdot \nabla) \, v)}{2}
\, \leq \, C \, (|\alpha|+\hnorme{I^s \widetilde{\nu}}{\sigma}) \,
(|\alpha|+\norme{\widetilde{w}}{2})
+ C \, (|\alpha|+\norme{\widetilde{\nu}}{\infty}) \, 
(|\alpha|+\norme{I^s \widetilde{w}}{2}) \ .
\end{displaymath}

Therefore, gathering these inequalities,
\begin{equation}
\begin{array}{lcl}
|\int_{\R 2}\,I^s \widetilde{w} &\times& I^s\,\textnormal{div}(\,b\,
\nabla^{\perp} \Pi)\,| \\[1ex]
 &\leq&C\,\norme{I^s \widetilde{w}}{2}^2\,
\norme{b_0}{\infty}\,(|\alpha|^2+\norme{\widetilde{\nu}}{\infty}^2)\\[1ex]
 &+&C\,\norme{I^s \nabla \widetilde{w}}{2}^2\,
(\norme{b_0}{\infty}+\varepsilon) \\[1ex]
 &+&C\,\hnorme{b}{s+\sigma}^2\,\big(\,\frac{1}{\varepsilon}+ 
\norme{b_0}{\infty}\,\big) \\[1ex]
 & &\times 
\big(\,|\alpha|^2+|\alpha|^2\,\norme{\widetilde{\nu}}{\infty}^2+|\alpha|^4+
(|\alpha|^2+\norme{\widetilde{\nu}}{\infty}^2)\,\norme{\widetilde{w}}{2}^2
+\norme{\nabla \widetilde{w}}{2}^2\,\big) \\[1ex]
 &+&C\,\norme{b_0}{\infty}\,\big(\,|\alpha|^2\,
(1+\norme{b_0}{\infty}^2+|\alpha|^2+\norme{\widetilde{\nu}}{\infty}^2) \\[1ex]
 & &+(|\alpha|^2+\hnorme{I^s \widetilde{\nu}}{\sigma}^2) \,
     (|\alpha|^2+\norme{\widetilde{w}}{2}^2)
    +\hnorme{b}{s+\sigma}^2(|\alpha|^2+\norme{\nabla \widetilde{w}}{2}^2)\, \big)\ .
\end{array}
\end{equation}
\end{itemize}

Putting all these points together and integrating yields
\begin{equation}
\begin{array}{lll}
&\norme{I^s \widetilde{w}}{2}^2 &
\\[1ex]
+& C \, \int_0^t \norme{I^s \nabla \widetilde{w}}{2}^2 &
\big( \, 1-\norme{b_0}{\infty} \, (1+|\alpha|)
- \sup_{[0,t]} \norme{I^s \widetilde{\nu}}{2} \, \big)
\\[1ex]
\leq & \norme{I^s \widetilde{w}_0}{2}^2 &
\\[1ex]
+& C \, \int_0^t \norme{I^s \widetilde{w}}{2}^2
&\big(1+\norme{b_0}{\infty} \int_0^t \norme{\widetilde{\nu}}{\infty}^2+
|\alpha|\,(1+\norme{b_0}{\infty} \int_0^t \norme{\widetilde{\nu}}{\infty}^2
|\alpha|\,\norme{b_0}{\infty})\big)
\\[1ex]
+& C \, \int_0^t \norme{\widetilde{\nu}}{\infty}^2 
&\big( \, \norme{b_0}{\infty}|\alpha|^2 \\[1ex]
& & +(|\alpha|^2 + \sup_{[0,t]} \norme{\widetilde{w}}{2}^2) \,
(1+\norme{b_0}{\infty}) \, \sup_{[0,t]}\hnorme{b}{s+\sigma}^2 \,\big)
\\[1ex]
+& C \, \int_0^t \norme{\nabla \widetilde{w}}{2}^2
&\big(\,|\alpha|+(1+\norme{b_0}{\infty}) \, 
\sup_{[0,t]} \hnorme{b}{s+\sigma}^2 \, \big)
\\[1ex]
+& C \, \int_0^t \gnorme{\widetilde{w}}{2}^2
&\big(\,|\alpha|+|\alpha|^2 \, \norme{b_0}{\infty}
+\sup_{[0,t]} \norme{I^s \widetilde{\nu}}{2} \\[1ex]
& & +|\alpha|^2 \, (1+\norme{b_0}{\infty}) \, \sup_{[0,t]} \hnorme{b}{s+\sigma}^2 \, \big)
\\[1ex]
+& C \, \int_0^t \hnorme{I^s \widetilde{\nu}}{\sigma}^2
&\norme{b_0}{\infty} \, (|\alpha|^2+\sup_{[0,t]} \norme{\widetilde{w}}{2}^2)
\\[1ex]
+& C \quad t &
|\alpha| \, (1+|\alpha|^3) \, (1+\norme{b_0}{\infty}) \,
(1+\sup_{[0,t]} \hnorme{b}{s+\sigma}^2)
\ .
\end{array}
\end{equation}

A Gronwall-type argument achieves the proof.
\end{proof}

\subsection{Estimate for convergence}

We now establish an estimate on the difference of two solutions of equations of type \eqref{vorticitelin}, intending to prove convergence of our iterative scheme and uniqueness of solutions of \eqref{complete}.

For $i=1,2$, consider
\begin{equation} \label{vorticitelindelta}
\partial_t \, \widetilde{w}_i 
      - (\mathcal{L} - \alpha \, \Lambda ) \, \widetilde{w}_i 
      + (\, \widetilde{\nu}_i \cdot \nabla \, ) \, \widetilde{w}_i 
=  \textnormal{div} \, \big( b_i \, (\nabla w_i
                                          + \nabla^{\perp} \Pi_i ) \big)
\end{equation}
where $\mathcal{L}$ and $\Lambda$ are as in \eqref{LandLambda}, $b_i$, 
$\widetilde{\Omega}_i$ are real functions, $\alpha \in \R{}$,
\begin{displaymath}
\begin{array}{rclcrcl}
\widetilde{v}_i & = & K_{BS} \star \widetilde{w}_i & , &
\widetilde{\nu}_i & = & K_{BS} \star \widetilde{\Omega}_i \ , \\
v_i & = & \alpha \, v^G + \widetilde{v}_i & , &
\nu_i & = & \alpha \, v^G + \widetilde{\nu}_i \ , \\
w_i & = & \alpha \, G + \widetilde{w}_i & , &
\Omega_i & = & \alpha \, G + \widetilde{\Omega}_i \ , 
\end{array}
\end{displaymath}
and $\nabla \Pi_i$ is obtained by solving
\begin{equation}
\textnormal{div} \, \big((1+b_i) \, \nabla \Pi_i \big)  
=\textnormal{div} \, \big((1 +b_i) \, \triangle v_i 
                         - (\nu_i \cdot \nabla) \, v_i \big) \ .
\end{equation}

Note that we choose to write $\widetilde{\nu}=K_{BS} \star \widetilde{\Omega}$ to stress the symmetry of the hypotheses on $\widetilde{\Omega}$ and 
$\widetilde{w}$.

For concision's sake, we note $\delta f=f_2-f_1$, for any functions $f_1$, $f_2$.

\begin{prop} $ $ \\
Let $\alpha \in \R {} $, $K > 0$, $\sigma > 2$, $0 < \eta < s <1$ and
$\max (\frac{2}{\eta},4) < p < +\infty$. \\
There exists $\varepsilon_0> 0$ 
and, for $K', T > 0$, there exists $C > 0$  such that if $\widetilde{w}_1$, 
$\widetilde{w}_2$ satisfy \eqref{vorticitelindelta} with
\begin{enumerate}
\item $\ \gnorme{b_0}{4} \leq \varepsilon_0 \, , \ 
\gnorme{b_0}{\infty} \leq \varepsilon_0$

\item  for $0 < t < T$, for $i=1,2$, for any $1 \leq r \leq +\infty$,
\begin{displaymath}
\begin{array}{lcrclcr}
\norme{b_i (t)}{r} &\leq & \norme{b_0}{r} \, e^{-\frac tr} &,&
\gnorme{b_i (t)}{r} &\leq & K \, \gnorme{b_0}{r} \, e^{-\frac tr}
\end{array}
\end{displaymath}

\item for $0 < t < T$, for $i=1,2$, $\ \hnorme{b_i (t)}{\sigma} \leq K'$

\item for $0 < t <T$, for $i=1,2$,
\begin{displaymath}
\begin{array}{rcl}
\gnorme{\widetilde{\Omega}_i (t)}{2}^2 + 
\int_0^t \gnorme{\nabla \widetilde{\Omega}_i}{2}^2 &\leq& \varepsilon_0\\[1ex]
\hnorme{\widetilde{\Omega}_i (t)}{s} +
\int_0^t \hnorme{\nabla \widetilde{\Omega}_i}{s}^2 &\leq& K'
\end{array}
\end{displaymath}

\item for $0 < t <T$, for $i=1,2$,
\begin{displaymath}
\begin{array}{rcl}
\gnorme{\widetilde{w}_i (t)}{2}^2 + 
\int_0^t \gnorme{\nabla \widetilde{w}_i}{2}^2 &\leq& \varepsilon_0\\
\hnorme{\widetilde{w}_i (t)}{s} + 
\int_0^t \hnorme{\nabla \widetilde{w}_i}{s}^2 &\leq& K'
\end{array}
\end{displaymath}
\end{enumerate}
then for $0 < t < T$,
\begin{eqnarray}\label{vorticitelinconvergence}
\gnorme{\delta \widetilde{w} (t)}{2}^2 &+& 
C \, \int_0^t \gnorme{\delta \widetilde{w}}{2}^2 
+ \gnorme{\nabla (\delta \widetilde{w})}{2}^2 
+ \gnorme{|\xi| \, (\delta \widetilde{w})}{2}^2 \nonumber \\
&\leq & C \, \int_0^t 
(1+ \gnorme{\widetilde{w}_1}{p}^2 + \hnorme{\nabla \widetilde{w}_1}{\eta}^2)
\  (\gnorme{\delta b}{p}^2 + \gnorme{\delta \widetilde{\Omega}}{2}^2) \ .
\end{eqnarray}
\end{prop}

\begin{proof}
Combining \eqref{vorticitedelta} for $i=1,2$, we derive
\begin{equation} \label{vorticitedelta}
\begin{array}{rcl}
\partial_{t} (\delta \widetilde{w}) 
&-&(\mathcal{L}-\alpha \, \Lambda) \, (\delta \widetilde{w})
+\big(\widetilde{\nu}_2 \cdot \nabla \big) \, (\delta \widetilde{w})
- \textnormal{div} \, \big( b_2 \, (\nabla (\delta \widetilde{w})
                                 + \nabla^{\perp} \Pi ) \big) \\[1ex]
&=& -((\delta \widetilde{\nu})\cdot \nabla) \, \widetilde{w}_1
    + \textnormal{div} \, \big( (\delta b) \nabla w_1 \big) \\[1ex]
& & + \, \textnormal{div} \, (b_2 \nabla^{\perp} R)
    + \, \textnormal{div} \, (b_2 \nabla^{\perp} (\delta S))
    + \, \textnormal{div} \, \big( (\delta b) \nabla^{\perp} \Pi_1 \big)
\end{array}
\end{equation}
with $\Pi$, $R$, $S_1$ and $S_2$ obtained by solving 
\begin{displaymath} \left\{
\begin{array}{lcl}
\textnormal{div} \, \big((1+b_2) \nabla \Pi  \big)
&=&\textnormal{div} \, 
\big((1+b_2) \triangle (\delta \widetilde{v})
-(\widetilde{\nu}_2 \cdot \nabla) \, (\delta \widetilde{\nu})\big)
\\
\textnormal{div} \, \big((1+b_2) \nabla R \big)
&=&\textnormal{div} \, 
\big(-\alpha (v^G \cdot \nabla) \, (\delta \widetilde{v})
+(\delta b) \triangle \widetilde{v}_1
-((\delta \widetilde{\nu}) \cdot \nabla) \, \widetilde{v}_1\big)
\\
\textnormal{div} \, \big((1+b_i) \nabla S_i \big)
&=& \textnormal{div} \, 
\big((1+b_1)\triangle v_1 - (\nu_1 \cdot \nabla) \, v_1 \big) \  ,
\qquad \textnormal{for} \ i=1,2 \ .
\end{array} \right.
\end{displaymath}
Note that $\Pi_2=\Pi+R+S_2$ and $\Pi_1=S_1$.

Our strategy is again to multiply \eqref{vorticitedelta} by $G^{-1} (\delta \widetilde{w})$, integrate on $\R 2$ and estimate each term arising to bound $\frac{d}{dt} \gnorme{\delta \widetilde{w}}{2}^2$. We deal with each term coming from the left member of equality \eqref{vorticitedelta} as we did in the linearized vorticity equation \eqref{vorticitelin}. Let us only show how to deal with the other terms.

First of all, let us emphasize that $\norme{G^r\,\nabla (G^{r'}\,f)}{2}^2$ 
is controled by $\gnorme{\nabla f}{2}^2 + \gnorme{|\xi| \, f}{2}^2$ provided that $r + r' = - \frac12$.

\begin{itemize}
\item First integrating by parts and applying H\"older's inequalities, we obtain
\begin{displaymath}
|\int_{\R 2} G^{-1} (\delta \widetilde{w}) \ 
\textnormal{div} (\widetilde{w}_1 (\delta \widetilde{\nu}))|
\leq \norme{G^{\frac12} \, \nabla (G^{-1} \, (\delta \widetilde{w}))}{2}
\,\norme{G^{-\frac12} \widetilde{w}_1}{p} \norme{\delta \widetilde{\nu}}{q}
\end{displaymath}
where $2 < q < + \infty$ is such that $\frac 1p + \frac 1q = \frac 12$.
Then using \eqref{BiotSavartLp}
\begin{equation}
|\int_{\R 2} G^{-1} (\delta \widetilde{w}) \ 
\textnormal{div} (\widetilde{w}_1 (\delta \widetilde{\nu}))|
\,\leq\,\varepsilon\,
\norme{G^{\frac12}\,\nabla (G^{-1}\,(\delta \widetilde{w}))}{2}^2
+\frac{C}{\varepsilon}\,\gnorme{\widetilde{w}_1}{p}^2
            \ \gnorme{\delta \widetilde{\Omega}}{2}^2
\end{equation}
where $\varepsilon$ is intended to be chosen small enough.

\item Similarly, with the same $q$, we also have
\begin{displaymath}
|\int_{\R 2} G^{-1}(\delta \widetilde{w}) \ 
\textnormal{div} ((\delta b) \nabla w_1)|
\leq \norme{G^{\frac12} \, \nabla (G^{-1}\,(\delta \widetilde{w}))}{2}
\, \norme{G^{-\frac12} (\delta b)}{p} \norme{\nabla w_1}{q}
\end{displaymath}
thus using Sobolev's embeddings
\begin{equation}
|\int_{\R 2} G^{-1}(\delta \widetilde{w}) \ 
\textnormal{div} ((\delta b) \nabla w_1)|
\leq \varepsilon \,
\norme{G^{\frac12}\,\nabla(G^{-1}\,(\delta \widetilde{w}))}{2}^2
+\frac{C}{\varepsilon}\,\hnorme{\nabla w_1}{\eta}^2
                      \ \gnorme{\delta b}{p}^2
\end{equation}
where $\varepsilon$ is again intended to be chosen small enough.

\item In quite the same way, we obtain
\begin{displaymath}
|\int_{\R 2} G^{-1} (\delta \widetilde{w}) \ 
\textnormal{div} (b_2 \nabla^{\perp} R)| \,
\leq\,C\,\norme{G^{\frac12}\,\nabla (G^{-1}\,(\delta \widetilde{w}))}{2}
\gnorme{b_0}{\infty} \norme{\nabla R}{2}\ .
\end{displaymath}
Thus with pressure estimate \eqref{pressionLp}, estimate on Biot-Savart law \eqref{BiotSavartLp}, H\"older's inequalities and Sobolev's embeddings we can derive
\begin{eqnarray}
|\int_{\R 2} G^{-1} (\delta \widetilde{w}) \ 
\textnormal{div} (b_2 \nabla^{\perp} R)|
& \leq & \varepsilon \,
\norme{G^{\frac12}\,\nabla (G^{-1}\,(\delta \widetilde{w}))}{2}^2
\nonumber \\
&  & + \, \frac{C}{\varepsilon}\,\gnorme{b_0}{\infty}^2
\hnorme{\nabla \widetilde{w}_1}{\eta}^2 \  \norme{\delta b}{p}^2
\nonumber \\
&  & + \, \frac{C}{\varepsilon}\,|\alpha| \gnorme{b_0}{\infty}^2 
\  \norme{\delta \widetilde{w}}{2}^2
\nonumber \\
&  & + \, \frac{C}{\varepsilon}\,\gnorme{b_0}{\infty}^2
\gnorme{\widetilde{w}_1}{2}^2 \ \gnorme{\delta \widetilde{\Omega}}{2}^2
\end{eqnarray}
where $\varepsilon$ is once again intended to be chosen small enough.

\item Using estimate \eqref{pressiondelta} instead of estimate \eqref{pressionLp} and inequalities \eqref{BiotSavartLp} and \eqref{BiotSavartGrad}, since $\frac 1p < \frac 12-\frac 1p$, we can obtain
\begin{eqnarray}
|\int_{\R 2} G^{-1} (\delta \widetilde{w}) \ 
\textnormal{div} (b_2 \nabla^{\perp} (\delta S))|
& \leq & \varepsilon \,
\norme{G^{\frac12}\,\nabla (G^{-1}\,(\delta \widetilde{w}))}{2}^2
\nonumber \\
&+&\! \frac{C}{\varepsilon}\,\gnorme{b_0}{\infty}^2\ 
\norme{\delta b}{p}^2
\nonumber \\
& &\! \! \! \! \! \times \big[
\hnorme{\nabla w_1}{\eta}^2 +\gnorme{\Omega_1}{2}^2 \gnorme{w_1}{p}^2 
\big]
\end{eqnarray}
where $\varepsilon$ is still intended to be chosen small enough.

\item At last, integrating by parts and applying H\"older's inequalities, we have
\begin{displaymath}
|\int_{\R 2} G^{-1} (\delta \widetilde{w}) \  
\textnormal{div}((\delta b) \nabla^{\perp} \Pi_1|
\leq \norme{G^{\frac12}\,\nabla (G^{-1} \, (\delta \widetilde{w}))}{2}
\, \gnorme{\delta b}{p} \norme{\nabla^{\perp} \Pi_1}{q} \ .
\end{displaymath}
where $2< q <+\infty$ is again such that $\frac 1p + \frac 1q = \frac 12$.
Again with pressure estimate \eqref{pressionLp}, estimates on Biot-Savart law \eqref{BiotSavartLp} and \eqref{BiotSavartGrad}, H\"older's inequalities and Sobolev's embeddings we can derive
\begin{eqnarray}
|\int_{\R 2} G^{-1} (\delta \widetilde{w}) \  
\textnormal{div}((\delta b) \nabla^{\perp} \Pi_1|
&\leq&\varepsilon \,
\norme{G^{\frac12}\,\nabla (G^{-1}\,(\delta \widetilde{w}))}{2}^2
\nonumber \\
&+&\! \frac{C}{\varepsilon}\,\gnorme{\delta b}{p}^2
\nonumber \\
& &\! \! \! \! \! \times \big[ 
\hnorme{\nabla w_1}{\eta}^2+\gnorme{\Omega}{2}^2 \gnorme{w_1}{p}^2 
\big]
\end{eqnarray}
where $\varepsilon$ is once again intended to be chosen small enough.
\end{itemize}

Gathering everything and integrating yield \eqref{vorticitelinconvergence}.
\end{proof}

\section{Main results}

We now use our various estimates to derive our main results.

\subsection{Existence and uniqueness}

Before proving a result of existence and uniqueness of solution of equations \eqref{complete}, we state a lemma that will make a link between norm estimates and convergence of the iterative scheme.

\begin{lemme}
Let $T >0$ and $p > 1$.
Let $(f_k)$ be a sequence in $L^{\infty}(0,T;\R +)$, and $(g_k)$ a bounded sequence in $L^p(0,T;\R+)$ be such that for $0<t<T$ and $k \in \N{}$,
\begin{displaymath}
f_{k+1}(t) \leq \int_0^t f_k\,g_k \ .
\end{displaymath}
Then $(f_k)$ is uniformally summable, namely, for $0<t<T$, 
\begin{displaymath}
\sum_{k \geq 0}f_k(t) \leq C_T \ .
\end{displaymath}
\end{lemme}

\begin{proof}
Using H\"older inequalities and iterating yield
$f_k(t)\, \leq \, K \, C^k \,(\frac{t^k}{k!})^{1-\frac 1p}$,
 where $C$ is a bound for $(g_k)$ in $L^p(0,T;\R+)$.
\end{proof}

We can now prove the existence and uniqueness parts of {\bf Theorem 1}.

\begin{proof}
\begin{itemize}
\item \emph{Existence.} We build a sequence 
$((b_k,\widetilde{w}_k))_{k \in \N*}$ such that, for any $k \in \N*$,
\begin{displaymath}
\left\{
\begin{array}{lcl}
\partial_t \, b_{k+1} + \big((v_k-\frac12 \xi)\cdot \nabla \big)\,b_{k+1}
&=&0 \\
\partial_t \widetilde{w}_{k+1} 
- (\mathcal{L}-\alpha\,\Lambda)\, \widetilde{w}_{k+1} 
+ \big(\widetilde{v}_k \cdot \nabla \big)\, \widetilde{w}_{k+1} 
&=&\textnormal{div}\, \big(b_k\,(\nabla w_{k+1}
+ \nabla^{\perp} \Pi_{k+1})\big)
\end{array}
\right.
\end{displaymath}
where $\mathcal{L}$ and $\Lambda$ are as in \eqref{LandLambda}, 
$(\widetilde{v}_k)$ is obtained from $(\widetilde{w}_k)$ by the Biot-Savart law, 
\begin{displaymath}
v_k=\alpha\,v^G+\widetilde{v}_k \ ,\quad
w_k=\alpha\,G  +\widetilde{w}_k \ ,\qquad \textrm{for} \quad k \in \N*,
\end{displaymath}
and $(\nabla \Pi_k)$ is obtained by solving, for any $k \in \N*$,
\begin{displaymath}
\textnormal{div} \,\big((1+b_k)\,\nabla \Pi_{k+1} \big)  
=\textnormal{div} \,\big((1+b_k)\,\triangle v_{k+1} 
                         -(v_k \cdot \nabla)\,v_{k+1}\big)
\end{displaymath}
with initial data $(b_0,\widetilde{w}_0)$. For $k=0$, we solve the system with $\widetilde{v}_k(t)\equiv 0$ and $b_k(t)\equiv 0$.

Let us show how we propagate bounds on $(b_k,\widetilde{w}_k)$.

{\bf Step 1} Fix $K_0>0$ and choose $\varepsilon_0 >0$ small enough. We can propagate \begin{enumerate}
\item for any $1 \leq p \leq +\infty$, $2\leq q \leq +\infty$, thanks to {\bf Proposition 5},
\begin{displaymath}
\begin{array}{lcrclcr}
\norme{b_k(t)}{p}  &\leq & \norme{b_0}{p} \, e^{-\frac tp} &,&
\gnorme{b_k(t)}{q} &\leq & \gnorme{b_0}{q} \,  e^{-\frac tq} \, e^{K_0} 
\end{array}
\end{displaymath}

\item and thanks to {\bf Prosposition 8},
\begin{eqnarray}
\gnorme{\widetilde{w}_k(t)}{2}^2
 &+& C_{K_0} \int_0^t \, \big( \gnorme{\widetilde{w}_k}{2}^2 
 + \gnorme{\nabla \widetilde{w}_k}{2}^2
 + \gnorme{|\xi| \widetilde{w}_k}{2}^2 \big)
\nonumber\\
 &\leq& C_{K_0} \,
(\gnorme{\widetilde{w}_0}{2}^2+\gnorme{b_0}{4})
\nonumber
\end{eqnarray}
\end{enumerate}
which provides us, thanks to {\bf Proposition 1} and Sobolev's embeddings,
\begin{displaymath}
\begin{array}{c}
\norme{\widetilde{v}_k}{8}\,\leq\,C\,\norme{\widetilde{w}_k}{\frac 85}\,
\leq\,C\,\gnorme{\widetilde{w}_k}{2}\,\leq\,K_0\\
\int_0^t \norme{\widetilde{v}_k}{\infty}^2\,
\leq\,C\int_0^t\,\big(\gnorme{\widetilde{w}_k}{2}^2 
+\norme{\nabla \widetilde{w}_k}{2}^2\big)\,
\leq\,\min(\frac{1}{24},K_0) \ .
\end{array}
\end{displaymath}

{\bf Step 2} Again choosing $\varepsilon_0$ small enough independently of $t$,and using {\bf Proposition 2}, we can obtain, when $0<s<1$ and $1<s''<2-s$,
\begin{displaymath}
\begin{array}{cclcl}
\int_0^t \hnorme{\nabla v_k}{1}&\leq&
C(t+t^{\frac 12}\,(\int_0^t \hnorme{\nabla \widetilde{v}_k}{1}^2)^{\frac 12})
&\leq&C(t+\int_0^t \hnorme{\widetilde{w}_k}{1}^2)\\
\norme{I^s \widetilde{v}_k(t)}{2}&\leq&C\ \gnorme{\widetilde{w}_k(t)}{2}
&\leq&C_{K_0}\,\varepsilon_0\\
\int_0^t \hnorme{I^s \widetilde{v}_k}{s''}^2&\leq&
C\,\int_0^t (\gnorme{\widetilde{w}_k}{2}^2+\norme{\nabla \widetilde{w}_k}{2}^2)
&\leq&K_0 \ 
\end{array}
\end{displaymath}
and propagate, for $0<t<T$,
\begin{enumerate}

\item when $1+s<s'<2$, thanks to {\bf Prosposition 6},
\begin{displaymath}
\hnorme{b_k(t)}{s'}\leq C_{K_0,T}
\end{displaymath}

\item and thanks to {\bf Prosposition 9},
\begin{displaymath}
\norme{I^s \widetilde{w}_k(t)}{2}^2 
+C \,\int_0^t \norme{I^s \nabla \widetilde{w}_k}{2}^2 \,
\leq\,C_{K_0,T}
\end{displaymath}
\end{enumerate}
which provides us, for $0<t<T$,
\begin{displaymath}
\int_0^t\,\hnorme{\nabla v_k}{s+1}\,\leq\,C\,(t+\int_0^t 
(\norme{\widetilde{w}_k}{2}^2+\norme{I^s \nabla \widetilde{w}_k}{2}^2))
\,\leq\,C_{K_0,T}
\end{displaymath}
thus, thanks to {\bf Prosposition 6}, for $0<t<T$,
\begin{displaymath}
\hnorme{b_k(t)}{s+2}\,\leq\,C_{K_0,T} \ .
\end{displaymath}

{\bf Step 3} We now prove the convergence of the scheme. Set $(\delta b)_k=b_{k+1}-b_k$ and $(\delta \widetilde{w})_k=\widetilde{w}_{k+1}-\widetilde{w}_k$. Choose $\max(4,\frac2s)<p<q$. {\bf Propositions 7 \& 10} give us, for $T >0$, for any $0<t<T$, for any $k \in \N*$,
\begin{equation} \label{convergence}
\begin{array}{l}
\gnorme{(\delta b)_{k+1} (t)}{p}^2 + 
\gnorme{(\delta \widetilde{w})_{k+1}(t)}{2}^2 \\[1ex]
\ \leq \ C_T \int_0^t \ 
(1+\gnorme{\widetilde{w}_k}{p}^2+\hnorme{\nabla \widetilde{w}_k}{\eta}^2)
\ \ (\gnorme{(\delta b)_k}{p}^2+\gnorme{(\delta \widetilde{w})_k}{2}^2)
\end{array}
\end{equation}
for some $0<\eta<s$ such that $\frac{2}{\eta}<p<+\infty$.

Now in order to apply {\bf Lemma 3} with $f_k=\gnorme{(\delta b)_k}{p}^2+\gnorme{(\delta \widetilde{w})_k}{2}^2$, remark that
\begin{itemize}
\item since $(G^{-\frac12} \widetilde{w}_k)$ is bounded in 
$L^{\infty}(\R+;L^2(\R 2))$ and $(\nabla (G^{-\frac12} \widetilde{w}_k))$ is bounded in $L^2(\R+; L^2(\R 2))$, \, $(G^{-\frac12} \widetilde{w}_k)$ is bounded in $L^r(\R+; L^p(\R 2))$, for some $2<r<+\infty$, by interpolation and Sobolev embeddings
\item since $(\widetilde{w}_k)$ is bounded in $L^{\infty}(\R+;L^2(\R 2)) \cap L^2(0,T;H^{s+1}(\R 2))$, \, $(\nabla \widetilde{w}_k)$ is bounded in 
$L^{r'}(0,T;H^{\eta}(\R 2))$ for some $2<r'<+\infty$, by interpolation.
\end{itemize}

Thus $(b_k)$ converges in $L^p_w(\R 2)$ and $(\widetilde{w}_k)$ in 
$L^2_w(\R 2)$, locally uniformally in $t$.

This implies at once that $(b_k)$ and $(\widetilde{w}_k)$ also converges in $L^2(\R 2)$. Now by interpolation
\begin{itemize}
\item since $(b_k)$ is bounded in $H^{s+2}(\R2)$, $(b_k)$ converges in 
$H^{\eta'}(\R2)$, for any $0<\eta'<s+2$
\item since $(\widetilde{w}_k)$ is bounded in $H^{s}(\R2)$, $(\widetilde{w}_k)$ converges in $H^{\eta''}(\R2)$, for any $0<\eta''<s$.
\end{itemize}
These properties enable us to take the limit in the sequence of equations.

Note that we recover the regularity on the limit by a mere application of Fatou's lemma.

\item \emph{Uniqueness.} We obtain a bound similar to \eqref{convergence} for the difference of two solutions. Then Gronwall lemma gives the result.
\end{itemize}
\end{proof}

\subsection{Asymptotic behaviour}

We now state the asymtotic part of {\bf Theorem 1}. Note that under the
hypotheses of {\bf Theorem 1}, the following asumptions are fulfilled.

\begin{theoreme}
Let $\alpha \in \R{}$. For any $0< \gamma <\frac12$, there exist 
$\varepsilon_0 > 0$ and $K, K' >0$ such that if $(b,\widetilde{w})$ is a 
solution of \eqref{complete} with initial data $(b_0,\widetilde{w}_0)$, such that
\begin{displaymath}
\gnorme{b_0}{2} \leq \varepsilon_0 \, , \quad 
\gnorme{b_0}{\infty} \leq \varepsilon_0 \, ,
\end{displaymath}
and for $t>0$, $\quad \gnorme{\widetilde{w}(t)}{2}^2
+\int_0^t \gnorme{\nabla \widetilde{w}}{2}^2\ \leq\ \varepsilon_0 \ ,$
\begin{displaymath}
\begin{array}{lclcrclc}
\gnorme{b(t)}{2} &\leq& K\,\gnorme{b_{0}}{2}\,e^{-\frac t2} &,&
\gnorme{b(t)}{\infty} &\leq& K\,\gnorme{b_0}{\infty} &,
\end{array}
\end{displaymath}
then, for $t>0$,
$\gnorme{\widetilde{w}(t)}{2} \leq K'\,e^{-\gamma\,t}\ 
(\gnorme{\widetilde{w}_0}{2}+\gnorme{b_0}{2})$.
\end{theoreme}

\begin{proof}
Let $0<\gamma<\gamma'<\frac12$. Let us recall \eqref{oscillateur}:
\begin{equation}
\int_{\R2} G^{-1}\,\widetilde{w}\ \mathcal{L}\,\widetilde{w} 
\ \leq\ -\gamma'\,\gnorme{\widetilde{w}}{2}^2
-(\frac12-\gamma')\,\big(\frac 13 \,\gnorme{\nabla \widetilde{w}}{2}^2 
+\frac 12 \,\gnorme{\frac{|\xi|}{4}\,\widetilde{w}}{2}^2\big) \ .
\end{equation}

Then we deal with the other terms of the vorticity equation as we did to obtain estimate \eqref{vorticitelinGLp}, except for the pressure term and
\begin{equation}
|\int_{\R 2} G^{-1}\,\widetilde{w}\ 
\widetilde{v} \cdot \nabla \widetilde{w}| \,
\leq \,C\,\gnorme{\widetilde{w}}{2}\,
(\gnorme{\widetilde{w}}{2}^2+\gnorme{\nabla \widetilde{w}}{2}^2) \ 
\end{equation}
obtained thanks to estimate \eqref{BiotSavartL8} and Sobolev's embeddings.

We treat the pressure term as follows
\begin{displaymath}
|\int_{\R 2} G^{-1}\,\widetilde{w}\ \textnormal{div}\,(b\,\nabla^{\perp} \Pi)|
\ \leq\ C\,\norme{G^{\frac 12}\,\nabla (G^{-1}\,\widetilde{w})}{2}\,
\norme{G^{-\frac 12}\,b\,\nabla^{\perp} \Pi}{2}
\end{displaymath}
with
\begin{displaymath}
\begin{array}{lcl}
\norme{G^{-\frac 12}\,b\,\nabla^{\perp} \Pi}{2}
&\leq&C\,\big(\gnorme{b_0}{2}\,e^{-\frac t2} 
\norme{\alpha \,(1+b)\,\triangle v^G-\alpha^2\,(v^G \cdot \nabla)\,v^G}{\infty}
\\[1ex]
&+&\!\!\!\gnorme{b_0}{\infty}\, 
\norme{(1+b)\,\triangle \widetilde{v}-\alpha\,((v^G \cdot \nabla)\,\widetilde{v}+(\widetilde{v}\cdot \nabla)\,v^G)-(\widetilde{v}\cdot\nabla)\,\widetilde{v}}{2}
\big)
\end{array}
\end{displaymath}
and $\quad \norme{(1+b)\,\triangle \widetilde{v}}{2}
\ \leq\ C\,(1+\gnorme{b_0}{\infty})\,\norme{\nabla \widetilde{w}}{2}$
\begin{displaymath}
\begin{array}{lclcl}
\norme{(v^G \cdot \nabla)\,\widetilde{v}}{2}&\leq&C\,
\norme{v^G}{\infty}\,\norme{\widetilde{w}}{2}
& &\\[1ex]
\norme{(\widetilde{v}\cdot \nabla)\,v^G}{2}&\leq&C\,\norme{\nabla v^G}{4}\,
\norme{\widetilde{v}}{4}
&\leq&C\,\gnorme{\widetilde{w}}{2}\\[1ex]
\norme{(\widetilde{v}\cdot\nabla)\,\widetilde{v}}{2}&\leq&C\,
\norme{\widetilde{v}}{\infty}\,\norme{\nabla \widetilde{v}}{2}&\leq&C\,
\big(\gnorme{\widetilde{w}}{2}+\norme{\nabla \widetilde{w}}{2}\big)\,
\norme{\widetilde{w}}{2} \ .
\end{array}
\end{displaymath}

This yields
\begin{equation}
\frac12 \frac{\textrm{d}}{\textrm{dt}}
\big(\gnorme{\widetilde{w}}{2}^2 \big)
+ \gamma \gnorme{\widetilde{w}}{2}^2 \,
+ \,C \,(\gnorme{\nabla \widetilde{w}}{2}^2
        +\gnorme{|\xi|\,\widetilde{w}}{2}^2)
\,\leq\,C\,\gnorme{b_0}{2}^2\,e^{-t}
\end{equation}
which integrating gives our result, since $2\,\gamma<1$.
\end{proof}

\bibliographystyle{plain}
\bibliography{Ref}

\end{document}